\documentclass[a4paper,12pt]{article}
\usepackage{amsmath}
\usepackage{bm}
\usepackage{makecell}
\usepackage{mathtools}
\usepackage{amssymb, amsbsy, amstext}
\usepackage{srcltx}
\usepackage{amscd,amsfonts,color}
\usepackage[figuresright]{rotating}
\usepackage{url}

\input amssym.def
\input amssym.tex

\newtheorem{lem}{Lemma}[section]%
\newtheorem{theorem}[lem]{Theorem}%
\newtheorem{prop}[lem]{Proposition}%
\newtheorem{rem}[lem]{Remark}%
\newtheorem{conj}[lem]{Conjecture}%

\def\a{\alpha} \def\b{\beta}  \def\d{\delta} 
    
\def\th{\theta}  
  
 \def\O{\Omega}

\def\v{{\text{\textasciicircum}}}

\def\di{\bigm|} \def\lg{\langle} \def\rg{\rangle}

\def\PSL{\hbox{\rm PSL}}  
\def\Aut{\hbox{\rm Aut\,}}  \def\Syl{\hbox{\rm Syl}}
 \def\soc{\hbox{\rm soc}} \def\Fix{\hbox{\rm Fix }}
  \def\mod{\hbox{\rm mod }}
  
 \def\PG{\hbox{\rm PG}} \def\PGL{\hbox{\rm PGL}}
  \def\GL{\hbox{\rm GL}}  \def\P\GL{\hbox{\rm P\GL}}
 \def\SL{\hbox{\rm SL}} \def\FF{{\hbox{\sf F\kern-.43emF}}}

\def\PSigmaU{{\rm P\Sigma U}}
\def\PGammaU{{\rm P\Gamma U}}
\def\Sym{\hbox{\rm Sym}}

\def\N{\hbox{\rm N}}
\def\C{\hbox{\rm C}}
\def\Z{\hbox{\rm Z}}

\def\soc{\hbox{\rm soc}}

\def\Sz{\hbox{\rm Sz}}
\def\Ree{\hbox{\rm Ree}}
\def\GU{\hbox{\rm GU}}
\def\PGU{\hbox{\rm PGU}}
\def\PSU{\hbox{\rm PSU}}
\def\PGammaU{{\rm P\Gamma U}}
\def\SO{\hbox{\rm SO}}
\def\SU{\hbox{\rm SU}}
\def\PSO{\hbox{\rm PSO}}
\def\PGO{\hbox{\rm PGO}}
\def\PSigmaU{{\rm P\Sigma U}}

\def\det{{\rm det}}

\def\ZZ{\mathbb{Z}}

\def\nd{\mathrel{\bigm|\kern-.7em/}} 
 \def\f{\noindent}
\def\qed{\hfill $\Box$} \def\demo{\f {\bf Proof}\hskip10pt}

\usepackage{titlesec}

\titleformat{\paragraph}
{\normalfont\normalsize\bfseries}{\theparagraph}{1em}{}
\titlespacing*{\paragraph}
{0pt}{3.25ex plus 1ex minus .2ex}{1.5ex plus .2ex}

\begin{document}

\begin{center} {\bf\large The Burness-Giudici Conjecture on
  Some Primitive \\ Groups with Socle $\PSU(3,q)$  }
\vskip 3mm
{\sc Huye Chen}\\
{\footnotesize
School of Mathematics, Guangxi University,   Nanning 530004, P. R. China}\\

{\sc Shaofei Du}\\
{\footnotesize
School of Mathematical Sciences, Capital Normal University, Beijing 100048, P. R. China}\\

{\sc Weicong Li}\\
{\footnotesize
School of Sciences, Great Bay University, Dongguan 523000, P. R. China}\\

\end{center}

\renewcommand{\thefootnote}{\empty}
 \footnotetext{{\bf Keywords} base of permutation group, Saxl graph, simple group}
  \footnotetext{E-mail addresses:  chenhy280@gxu.edu.cn (H. Chen), dushf@mail.cnu.edu.cn (S. Du),
  liweicong@gbu.edu.cn (W. Li).}

\begin{abstract}
Let $G$ be a transitive permutation group on $\O$ with two points $\a, \b\in\O$ such that
$G_\a\cap G_\b=1$. The  Saxl graph $\Sigma(G)$ of the pair $(G, \O)$ is the graph with  vertex set $\O$, while two vertices $\a', \b'$ are adjacent if and only if $G_{\a'}\cap G_{\b'}=1$.
It was conjectured by Burness and Giudici that
 the  Saxl graph $\Sigma(G)$ of  any primitive permutation group $G$
has the property that any two vertices have a common neighbor.
We focused on proving the conjecture for
  all  primitive  groups $G$ whose socle is a simple group of  Lie-type of rank $1$, that is, those with $\soc(G)\in \{\PSL(2,q),\PSU(3,q), \Ree(q),\Sz(q)\}$.
The case of $\soc(G)=\PSL(2,q)$ has been published in two papers. This paper will address most cases where $\soc(G)=\PSU(3,q)$, with the exception of a particularly intricate configuration in which the point stabilizer contains $\PSO(3,q)$. That specific configuration has been treated in a separate paper.
\end{abstract}

\section{Introduction}
Let $G\leq \Sym(\Omega)$ be a permutation group. Then a subset $\Delta$ of $\O$ is called a {\em base} if its pointwise stabilizer therein is trivial.
 The {\em base size} of $G$ on $\O$, denoted by $b(G)$, is the minimal cardinality of a base for $G$. A {\em base size set} $\Delta$ for a
 finite permutation group $G\leq \Sym(\Omega)$,  is a base for $G$ on the set $\O$ with $|\Delta|=b(G)$.

 Recently, Burness and Giudici \cite{BG} introduced a graph, called the {\em Saxl graph} which is related  to  bases of a permutation
group with $b(G)=2$. Let $G\leq \Sym(\Omega)$ be a permutation group on $\Omega$ with $b(G)=2$. The vertex set of a Saxl graph $\Sigma(G)$ (denoted by $\Sigma$, simply) of $G$ on $\Omega$ is just
$\Omega$ and two vertices are adjacent if and only if they form  a base size set.

 Henceforth, assume that $G$ acts transitively on $\Omega$. Fix a point $\alpha \in \Omega$. The orbits of the stabilizer $G_\alpha$ are called the {\em suborbits} of $G$ with respect to $\alpha$, where $\{\a\}$ is said  to be trivial.
  A pair ${\alpha, \beta}$ is a base for $G$ if and only if $G_{\alpha}$ acts regularly on the suborbit containing $\beta$, given that $b(G)=2$. Note that the Saxl graph $\Sigma(G)$ is vertex-transitive whenever $G$ is a transitive permutation group on $\Omega$. Moreover, the neighborhood $\Sigma_1(\alpha)$ of $\alpha$ in $\Sigma(G)$ is the union $\Gamma$ of the regular suborbits of $G$ relative to $\alpha$.
In their paper \cite{BG}, the authors investigated the valency, connectivity, hamiltonicity, and independence number of $\Sigma(G)$. They also proposed a conjecture, which can be stated as follows.

 \begin{conj}$($BG-Conjecture$)$ \label{BG1}
If $G$ is a finite primitive permutation group with $b(G)=2$, then every pair of vertices in its Saxl graph $\Sigma(G)$ has a common neighbor.
\end{conj}

Burness and Giudici \cite{BG} employed probabilistic methods to verify the conjecture for a broad class of almost simple groups. Their work specifically establishes the conjecture for symmetric and alternating groups
$G=S_n$ or $A_n$ (with $n>12$) when the point stabilizer is primitive. Their results also apply to certain primitive groups of diagonal type and to twisted wreath products of sufficiently large order; these cases have recently been extended by Huang (see Sections 5.8 and 5.9 of \cite{HPHD}).
Moreover, by using the Magma computational algebra system, they confirmed the conjecture for every primitive group of degree at most $4095$. The same paper \cite{BG} also verifies the conjecture for many of the sporadic simple groups.
Subsequently, Burness and Huang \cite{BH} proved the conjecture for almost simple primitive groups whose point stabilizers are soluble. In a separate article \cite{BH2}, the same authors further investigated the conjecture for primitive groups of product type.
Most recently, Lee and Popiel \cite{LP} provided a proof for the majority of affine-type groups with sporadic point stabilizers. The framework of the conjecture has since been expanded with the introduction of generalized Saxl graphs and Saxl hypergraphs in \cite{FH, LP1}.

To date, this conjecture is still open for simple groups of Lie type. Definitely, this is a huge and difficult  topic, especially considering that a full classification of base-two primitive almost simple groups remains out of reach.
 The first candidate  may be  those groups of rank one, that is $\soc(G)\in \{\PSL(2,q), \PSU(3,q), \Ree(q), \Sz(q)\}$.   The conjecture has been settled for groups with socle $\PSL(2,q)$ in \cite{BH, Chen-Du}. Accordingly, we focus on the remaining three families.
  We shall finish them with three papers.
The first, \cite{Chen-Du-Li}, treated the case where $\soc(G) = \PSU(3,q)$ and a point stabilizer contains $\SO(3,q)$.
The present paper establishes the conjecture for primitive permutation representations $G$ with $\soc(G) = \PSU(3,q)$ whose point stabilizer does not contain $\PSO(3,q)$ (see Theorem~\ref{main}).
Finally, the case of groups $G$ with socle $\Ree(q)$ or $\Sz(q)$ will be addressed in a separate work \cite{Chen-Du-2}.
Now we are ready to  state  the main result of this paper.
\begin{theorem} \label{main}  Let $G$ be a primitive group with socle $\PSU(3,q)$, whose
  point-stabilizer does not contain $\PSO(3,q)$.
 Then the BG-Conjecture holds for $G$.
  \end{theorem}

This paper is organized as follows. Following this introduction, Section 2 provides the necessary notation, basic definitions, and useful facts. The proof of Theorem \ref{main} is presented in Section 3. With the exception of certain cases that are settled within Sections 3 and 7, the proof reduces to studying a specific family of groups. This family will be examined in detail in Sections 4-6: general preparations are provided in Section 4, and the two main subcases are analyzed in Sections 5 and 6, respectively.

\section{Preliminaries}
Let $\FF_{q}$ be a finite field of order   $q=p^m$, where $p$ is a prime and $m$ is a positive integer. Let $\FF_{q^2}^*=\lg \xi\rg$  and $\FF_{q}^*=\lg \th\rg $ be the cyclic multiplicative groups.
   The Frobenius automorphism $f \in \operatorname{Aut}(\FF_{q^2})$, defined as $f(\xi) = \xi^p$, has order $2m$.

  Let $V$ be the $3$-dimensional unitary space of row vectors over $\FF_{q^2}$, equipped with a Hermitian form $(\cdot, \cdot)$. We denote by $\GU(3,q)$, $\SU(3,q)$,  $\PGU(3,q)$ and $\PSU(3,q)$ the general unitary, special unitary, projective general unitary, and projective special unitary groups, respectively.
  Note that $\PSU(3,q)$ is simple for $q > 2$. We recall some known results for $T := \PSU(3,q)$, whose order is given by $|T| = \frac{1}{d} (q^3+1)q^3(q^2-1)$ with $d:=(3,q+1)$.
  Throughout this paper, we let $\Syl_p(G)$ denote the set of sylow $p$-subgroups of a group $G$. Moreover, for any $a, b \in G$, we define the commutator as $[a, b] := a b a^{-1} b^{-1}$.
\vskip 3mm
(1) We choose a unitary basis ${w_1, w_2, w_3}$ whose Gram matrix is the antidiagonal matrix $J:=]1,1,1[$ (see \cite{Huppert}).
Using the automorphism $a\mapsto \overline{a}:=a^q$ $(a\in\FF_{q^2})$, we define $\GU(3,q):=\{g\in\GL(3,q^2)\mid gJ\overline{g}^T=J\}$. Let ${\v}X:=X\Z/\Z$ be the image of a subgroup $X$ of $\GU(3,q)$, where $\Z:=\Z(\GU(3,q))$ and use $\hat{g}$ to denote the image of an element $g\in \GU(3,q)$. In $\SU(3,q)$, define
\vspace{-5pt}$${\small\begin{array}{lll}
q(a,b)={{\left(
\begin{array}{ccc}
1&a&b \\
0&1&-\overline{a} \\
0&0&1\\
\end{array}
\right)  }},&  h(k)={{ \left(
\begin{array}{ccc}
k^{-q}&0&0 \\
0&k^{q-1}&0\\
0&0&k\\
\end{array}
\right)  }} &\, \tau={{ \left(
\begin{array}{ccc}
0&0&1 \\
0&-1&0\\
1&0&0\\
\end{array}
\right)  }},
\end{array}}$$
$$Q:=\{ \hat{q}(a,b)\mid a,b\in \FF_{q^2}, b+\overline{b}+a\overline{a}=0 \}, \quad H:=\{\hat{h}(k)\mid k\in\FF_{q^2}^*\}.$$
Thus, $Q$ is a Sylow $p$-subgroup of $T=\PSU(3,q)$ of order $q^3$; and $H$ is a cyclic subgroup of $T$ of order $\frac{q^2-1}{d}$,
with relations
 \vspace{-5pt}\begin{equation}\label{compu}
 \begin{array}{lll}
 &\hat{q}(a,b)\hat{q}(a',b')=\hat{q}(a+a',b+b'-a\overline{a'}), \quad
&\hat{\tau} \hat{h}(k)\hat{\tau}=\hat{h}(k^{-q}),\\
&\hat{h}(k)^{-1}\hat{q}(a,b)\hat{h}(k)=\hat{q}(k^{2q-1}a,k^{q+1}b),
&\hat{q}(a,b)^{-1}=\hat{q}(-a,\overline{b}).\\
\end{array}
\end{equation}
Further, the center of $Q$, given by $\Z(Q)=:Q_0=\{\hat{q}(0,b)\mid b\in\FF_{q^2}, b+\overline{b}=0 \}$, has order $q$.
 \vskip 3mm
 (2)  Let $\delta = [\mu, 1, 1] Z$, where $\Z := Z(\GU(3, q))=\{[\mu,\mu,\mu]\mid \mu\in\FF_{q^2}, \mu^{q+1}=1\}$; here $[\mu, 1, 1]$ denotes the diagonal matrix with entries $(\mu, 1, 1)$, and $|\delta|=d$. Hence, $\PGU(3,q)=T\rtimes \lg \d\rg $.
 The field automorphism of $T$ induced by $f$ is redenoted by $f$ again so that
 $\Aut(T)=T\rtimes \lg \d, f\rg .$
 \begin{prop}\label{fieldconjugate}{\rm \cite[p.77,79]{BGbook}}.
Let $T=\PSU(3, q)$ and $f$ the field automorphism. Then
\begin{enumerate}
\item[\rm(i)] Let $x\in\lg f\rg$ such that  $x$ is of order a prime $r$.  Then for any $g\in \PGU(3,q)$,
   $|gx|=r$ if and only if $gx=x^{g_1}$  for some $g_1\in \PGU(3,q)$.
\item[\rm(ii)] Let $z\in\PGU(3,q)$ be an element of prime order $r\ne p$.  Then $z^T=z^{\PGU(3,q)}$.
\end{enumerate}
\end{prop}

\begin{prop}\label{man}
{\rm \cite{Man}} \
Let $G$ be a transitive group on $\O$ and let $H=G_\a$
for some $\a\in \O$.
Suppose that $K\le G$ and at least one $G$-conjugate of
$K$ is contained in $H$. Suppose further that the set of
$G$-conjugates of $K$ which are contained in $H$ form $t$
conjugacy classes of $H$ with representatives $K_1$, $K_2$, $\cdots,$ $K_t$.
Then $K$ fixes $\sum_{i=1}^{t}|\N_G(K_i):\N_H(K_i)|$ points of $\O$.
\end{prop}

Let $fix(x,G/M)$ denote the number of fixed points of $x$ on $\O:=[G:M]$ and let
\vspace{-5pt}$$fpr(x,G/M)=\frac{fix(x,G/M)}{|\O|}$$
be the fixed point ratio of $x$.
The following proposition is modified from \cite{BG}.

\begin{prop}\label{prob}
Consider the transitive permutation representation of $G$ on $\O:=[G:M]$. Let $\check{Q}(G):=\frac{|M|}{|\O|}(\sum_{i=1}^k(|\lg x_i\rg|-1)\frac{\Fix(\lg x_i\rg)}{|\N_M(\lg x_i\rg )|})$, where
 $\mathcal{P}^*(M):=\{\lg x_1\rg,\lg x_2\rg,\cdots,\lg x_k\rg\}$ is the set of representatives  of conjugacy classes of subgroups of prime order in $M$ and $\Fix(\lg x_i\rg)$ is the number of fixed
 points of $\lg x_i\rg$ on $\O$.  Then any two vertices in $\Sigma(G)$ have a common neighbour, provided  $\check{Q}(G)< \frac 12$.
\end{prop}

\begin{rem}\label{big} Let $G$ and $G_0$ be primitive groups on $\O$ with $b(G)=b(G_0)=2$, where $G_0$ is a subgroup of $G$ with a point stabilizer contained in that of $G$. If the BG-Conjecture holds for $G$ acting on $\O$, then the BG-Conjecture also holds for $G_0$ acting on $\Omega$. The result follows from the inclusion $\Gamma(\a)\subseteq\Gamma_0(\a)$ of the respective sets of regular suborbits.
\end{rem}

\section{The proof of Theorem~\ref{main}}
 {\bf Proof}   Let $G$ be an almost simple group with socle $\soc(G)=T=\PSU(3,q)$, that is, $T\leq G\leq \Aut(T)$.
 Consider the primitive permutation representations of $G$ on the right cosets $[G:M]$ realtive to    maximal subgroups $\a=M$ of $G$, up to conjugacy. To prove Theorem \ref{main},
 we shall find a common  neighbor of $\a $ and  $\a^g$ for any $g\in G$ in the Saxl graph.

Suppose that  the group  $M$ is soluble. Then the  conjecture is true by  \cite[Theorem 1.1]{BH}. So in what follows, we assume that {\it $M$ is insolvable.}
 If $T\lesssim M$, then $G/M_G$ is cyclic, where $M_G=\bigcap\limits_{g\in G}M^g$. In this case,  $b(G)=1$, which is not our case. So
  we assume  that {\it $T\not\lesssim M$.}
   The case when $q=5$ may be vertified by   Magma directly and so assume that {\it $q\ge 7$} below.

Let  $M_0:=M\cap T$.  Then  by \cite[Table 8.5 and Table 8.6]{BGbook},  $M_0$ is a maximal subgroup of $T$ and
$$M_0\cong L\in \{ \SL(2,q)\rtimes \ZZ_{\frac{q+1}{d}},  \PSU(3,q').\ZZ_c,  \PSL(2,7), \PSL(2, 9),\PSO(3,q)\},$$ where $q=q'^e$, $c:=(\frac{q+1}{q'+1}, 3)$ and $e$ is an odd prime. Note that $M_0 \cong \PSL(2,7)$ for $q = p \equiv 3, 5, 6 \pmod{7}$ with $q \neq 5$; and $M_0 \cong \PSL(2,9)$ for $q = p \equiv 11, 14 \pmod{15}$. The first four cases will be dealt with separately.
\vskip 3mm
(1)   $M_0\cong\SL(2,q)\rtimes \ZZ_{\frac{q+1}{d}}$

\vskip 3mm Let  $G$ be an almost simple group with socle $\soc(G)=T=\PSU(3,q)$, $M$ a maximal subgroup of $G$ containing  $M_0=\SL(2,q)\rtimes \ZZ_{\frac{q+1}{d}}$.
In what follows, we shall show $b(G)\ne 2$.

 As a subgroup of $G$, the action of $T$ on  $\O=[G:M]$ may be identified with  its primitive  right multiplication   on  cosets  $\O_0=[T:M_0]$.
Clearly, there exist no regular suborbits of $G$ on $\O$, that is $b(G)\ne 2$, provided
 there exist no regular suborbits of $T$,  that is $b(T)\ne 2$.

 Now, the  action of $T$ on $\O_0=[T:M_0]$ is equivalent to its action on the set $\mathcal{P}_{non}$  of non-isotropic points in the unitary projective plane $\PG(2, q^2)$.
 We shall show $T_{\lg \a\rg} \cap T_{\lg \b\rg } \ne 1$ for any $\lg\a\rg,\lg\b\rg\in\mathcal{P}_{non}$ so that $b(T)\ne 2$. Since $T$ acts transitively on $\mathcal{P}_{non}$, we may fix $\alpha = (1, 0, 0)$ without loss of generality. For any $g \in \SU(3,q)$, the defining condition $g \overline{g}^T = E$ (where $E$ is the identity Gram matrix), when applied to the stabilizer $T_{\lg \a \rg}$, implies that
 $${\small T_{\lg \a \rg}=\lg \hat{g}\di  g=\left(
\begin{array}{cc}
a & 0 \\
0 & A \\
 \end{array} \right)  \in \SU(3,q),  A\in \GU(2,q), a|A|=1\rg .}$$
Now $\b$  is of form  $(0,0,1)$, $(0, 1, c_1)$ or $(1, b, c_2)$ for some $b,c_1,c_2\in\FF_{q^2}$ with $c_1^{1+q}+1\neq0$ and $1+b^{1+q}+c_2^{1+q}\neq0$.
Since ${\v}A$ acts nonregularly on $\lg(0,1)\rg$ $\lg (1, c_1)\rg $ and  $\lg (b, c_2)\rg $,  we know $T_{\lg \a \rg} \cap T_{\lg \b \rg}\ne 1$.

\vskip 3mm
(2) $M_0\cong\PSU(3,q').\ZZ_c$, which will be done in Sections 4-6, see Theorem~\ref{main1}.
\vskip 3mm
(3) $M_0\cong\PSL(2,7), \PSL(2,9)$, see Theorems~\ref{psl27} and ~\ref{psl29} in Section 7.

\vskip 3mm
 In summary, Theorem~\ref{main} is proved. \qed

\section{$M_0\cong\PSU(3,q').(\frac{q+1}{q'+1},3)$}
In this section, we mainly prove the following theorem.
\begin{theorem}\label{main1} Let $G$ be a primitive group with sloce $T=\PSU(3,q)$ acting on a set $\O$, such that a point stabilizer contains a subgroup  $M_0\cong\PSU(3,q').(\frac{q+1}{q'+1},3)$.
  Here, $q = p^m$ and $q' = p^{m'}$ with $m = m'e$, where $p$ is a prime and $e$ is an odd prime.
 Suppose that $b(G)=2$. Then
  the BG-Conjecture is true.
\end{theorem}
\demo Let $c=(\frac{q+1}{q'+1},3)$, $T = \PSU(3,q)$, and $G$ be a group such that $T \leq G \leq \PGammaU(3,q)$. Furthermore, let $M$ be a maximal subgroup of $G$ containing $M_0 = \PSU(3,q') \cdot \ZZ_{c}$.
 As shown in Table $8.5$ of \cite{Bray}, we have $T\leq G \le T :\lg \delta^c, f \rg$, where  $\delta$  is a diagonal automorphism of $T$ of order $d=(q+1,3)$.
By the definition of the unitary groups, with $q = p^m$ and $\gamma := f^m$, we have the conjugation relations $\delta^f = \delta^p$ and $\delta^\gamma = \delta^{-1}$.

Write $q' + 1 = 3k + t$, where $t \in \{0, 1, 2\}$ and $k$ is a nonnegative integer. Then  \vspace{-5pt}$$q+1=(3k+t-1)^e+1\equiv  (t-1)^e+1\equiv 0\pmod 3$$ if and only if $t=0$, and so  $(q'+1, 3)=(q+1, 3)$.

Since $\frac{q+1}{q'+1}=\frac{((q'+1)-1)^e+1}{q'+1}\equiv e \pmod{q'+1}$ and $e$ is an odd prime, we get  $c:=(\frac{q+1}{q'+1},3)=3$ if and only if $e=3$ and $(q'+1,3)=3$.

In the following two sections,  we deal with two cases $c=1$ and $e\geq 3$;  and $c=e=3$, separately. Therefore, Theorem~\ref{main1} is a direct consequence of Theorems~\ref{saxl3} and \ref{c=3proof}. \qed

\vskip 3mm
Combining Propositions 3.3.2, 3.3.3, 3.3.7, and Table B.4 in \cite{BGbook}, we obtain the following lemma.
\begin{lem}\label{PGU}
Let $T = \PSU(3,q)$ and $R =\PGU(3,q)=T\rtimes\lg \d\rg$, where $(3,q+1)=3$. For any element $x \in R$ of prime order $r$, the $R$-conjugacy class of $x$ is determined as follows:
\begin{enumerate}
    \item If $r \neq 3$, then $x$ is $R$-conjugate to an element of $T$.
    \item If $r = 3$, then $x$ is $R$-conjugate to one of the three diagonal matrices $z_A$, $z_B$, or $\delta=[\mu,1,1]\Z$, whose centralizers in $R$ are respectively
        $\C_R( z_{A}) \cong (\ZZ_{q+1} \times \ZZ_{q+1}) \rtimes \ZZ_3,
    \C_R(z_{B}) \cong \SL(2,q) \rtimes \ZZ_{q+1},$
     and $
    \C_R(\delta) \cong \ZZ_{q^2-q+1} \rtimes \ZZ_3.
    $
\end{enumerate}

\end{lem}

In the remaining of this section, some preparations are given. To prove Theorem~\ref{main1} via Proposition~\ref{prob}, we require information on the prime order subgroups of $\PGU(3,q)$.

\begin{lem}\label{primeproperties}
Let $T=\PSU(3,q)$ and $R=\PGU(3,q)$. Then every cyclic subgroup of prime order $r$ in $T$ is $T$-conjugate to $\lg z \rg$, where the canonical representation of $z$, $|\N_T(\lg z\rg )|$ and $|\N_R(\lg z\rg )|$
are  given in Table \ref{sy}.
\end{lem}
\demo Noting that $|T|=\frac 1d (q^3+1)q^3(q^2-1)$, for any prime dicisor $r$ of $|T|$ the possible cases for $r$ are $r\in \{2, 3, p\}$, or $r\di (q-1)$, $r\di (q+1)$ or $r\di (q^2-q+1)$. From Table~\ref{character}, which describes the conjugacy classes of elements in $T$ \cite{Simpson-Frame}, we obtain the following cases based on the value of $r$.
\vskip 3mm
(1) $r=2$ and $p$ is odd:
\vskip 3mm
Note that the centralizer of every involution is the stabilizer of a non-isotropic point. Since both $T$ and $R$ are transitive on the set of non-isotropic points, it follows that all subgroups isomorphic to $\ZZ_2$ are conjugate in the group. Moreover, for any such involution $z_2$, we have $|\N_T(\langle z_2 \rangle)| = q(q^2-1)(q+1)/d$ and $|\N_R(\langle z_2 \rangle)| = q(q^2-1)(q+1)$.

\vskip 3mm
(2) $r=p$:
\vskip 3mm
From Table \ref{character}, we get that $|\C_T(z_0)|=\frac{q^3(q+1)}{d}$ and there is only one conjugacy class subgroup $\lg z_0\rg$ of this form, which presented in Table \ref{sy}. Note that $\N_T(\lg z_0\rg)/\C_T(\lg z_0\rg)\lesssim \ZZ_{p-1}$. From Eq(\ref{compu}), we have $\hat{h}(k)^{-1}\hat{q}(0,1)\hat{h}(k)=\hat{q}(0,k^{q+1})$ and $(\hat{q}(0,1))^l=\hat{q}(0,l)$
for any $k\in\FF_{q^2}$ and positive integer $l$. Therefore, we have $|\N_T(\lg z_0\rg)|=\frac{q^3(q+1)(p-1)}{d}$. Furthermore, $|\C_R(\lg z_0\rg)|=q^3(q+1)$  followed from Proposition \ref{fieldconjugate}. This implies that $|\N_R(\lg z_0\rg)|=q^3(q+1)(p-1)$, as $\N_R(\lg z_0\rg)/\C_R(\lg z_0\rg)\lesssim \ZZ_{p-1}$.

Table \ref{character} shows that there are $d$ conjugacy classes of $z_l'$ in $T$, each with $|\C_T(\lg z_l' \rg)| = q^2$, and that these $d$ classes are fused into a single class under the action of $R$. This implies that $|\C_R(\lg z_l' \rg)| = q^2$.
Note that we may, without loss of generality, choose one of the elements
$z_l'$ to lie in $\PGO(3,q)\cong\PGL(2,q)$; denote this element by $z_0'$. Since $\ZZ_{p-1}\lesssim\N_{\PGO(3,q)}(z_0')$ and $\N_T(\lg z_0'\rg)/\C_T(\lg z_0'\rg)\lesssim \ZZ_{p-1}$,
we have $\N_T(\lg z_0'\rg)=q^2(p-1)$. This implies that $\N_R(\lg z_l'\rg)=q^2(p-1),$ for any $l$, since these elements are all conjugate in $R$. Consequently, we also have $\N_T(\lg z_l'\rg)=q^2(p-1)$ for any $l$.

\vskip 3mm

(3) $r=3$ and $3\mid q+1$:
\vskip 3mm
As can be seen from Table \ref{character}, $z_3^A$ exists only when $3 \mid q+1$, and $z_3^B$ exists only when $9 \mid q+1$, where $z_3^A$ and $z_3^B$ are of the form presented in Table \ref{sy}. And $\lg z_3^A\rg$ is not conjugate to $\lg z_3^B\rg$ in $\PSU(3,q)$.
Moreover, for $3\mid q+1$, there exist only these two conjugacy classes of elements of order $r=3$ in $\PSU(3,q)$.
 By Table \ref{character}, $|\C_T(z_3^A)|=(q+1)^2$ and $|\C_T(z_3^B)|=\frac{q(q+1)^2(q-1)}{3}$.
By Proposition \ref{fieldconjugate}, we get that $|T:\C_T(x)|=|R:\C_R(x)|$ for any $x\in T$ of order co-prime to $p$. Therefore, we get that $|\C_R(z_3^A)|=3(q+1)^2$ and $|\C_R(z_3^B)|=q(q+1)^2(q-1)$.

Note that $\N_T(\lg z_3^A\rg)/\C_T(\lg z_3^A\rg)\lesssim \ZZ_{2}$. This implies that $|\N_T(\lg z_3^A\rg)|=(q+1)^2s$, where $s=1,2$. Hence $(q+1)^2s\mid |\N_R(\lg z_3^A\rg)|$. From $|\C_R(z_3^A)|=3(q+1)^2$, we conclude  that $3(q+1)^2\mid |\N_R(\lg z_3^A\rg)|$. Therefore, $3s(q+1)^2\mid |\N_R(\lg z_3^A\rg)|$. Hence, the order must be $|\N_R(\lg z_3^A\rg)|=3s(q+1)^2$ with $s=1,2$.

By anlyzing the structure of the subgroup of $\PSU(3,q)$, we get that $\C_T(\lg z_3^B\rg)\cong \SL(2,q)\rtimes\ZZ_{\frac{q+1}{d}}$ and $\C_R(\lg z_3^B\rg)\cong \SL(2,q)\rtimes\ZZ_{q+1}$.
Consequently, the centralizer $\C_T(\lg z_3^B\rg)$ (resp. $\C_R(\lg z_3^B\rg)$)is a maximal subgroup of $T$ (resp. $R$). This implies that $|\N_T(\lg  z_3^B\rg)|=|\C_{T}(\lg z_3^B\rg)|=\frac{q(q+1)^2(q-1)}{3}$ and
$|\N_R(z_3^{B})|=|\C_R(z_3^{B})|=q(q+1)^2(q-1)$.

\vskip 3mm
(4) $r\di (q+1)$, $r\ne 2, 3$:
\vskip 3mm

By Table \ref{character}, for $r\mid q+1$ and $r\neq 2,3$, there exist only these two types of conjugate classes elements in $\PSU(3,q)$ of order $r$, with respective elements $z_{q+1}^A$ and $z_{q+1}^B$. For the type $z_{q+1}^B$, there is only one conjugate classes subgroup $\lg z_{q+1}^B\rg$ of this type. For the type $z_{q+1}^A$, there are at most $\frac{r-3}{2}$ (depending on the choices of $(1,a,b)$ such that $1+a+b\equiv 0(\mod r)$) conjugate classes subgroup $\lg z_{q+1}^A(a)\rg$ of this type, depending on the choices of $a$.
By Table \ref{character}, we have that
$|\C_T(z_{q+1}^A (a))|=\frac{(q+1)^2}{d}$ and
$|\C_T(z_{q+1}^B)|=\frac{q(q+1)^2(q-1)}{d}$.
By Proposition \ref{fieldconjugate}, we conclude that $|\C_R(\lg z_{q+1}^A(a)\rg)|=(q+1)^2$ and $|\C_R(\lg z_{q+1}^B\rg)|=q(q+1)^2(q-1)$. It follows that $|\N_R(\lg z_{q+1}^B\rg)|=q(q+1)^2(q-1)$, as a maximal subgroup of $R$ has order $q(q+1)^2(q-1)$.
And we can also get  that $|\N_T(\lg z_{q+1}^B\rg)|=\frac{q(q+1)^2(q-1)}{d}$, as a maximal subgroup of $T$ has order $\frac{q(q+1)^2(q-1)}{d}$.

By analyzing the structure of group $T$, we get that $\N_T(\lg z_{q+1}^A \rg)=\ZZ_{\frac{q+1}{d}}\times\ZZ_{q+1}\rtimes S$ where $S\lesssim S_3$ is a permutation matrix and $S$ is determined by the choices of $(1,a,b)$ where $1+a+b\equiv 0(\mod r)$.
So that $\N_R(\lg z_{q+1}^A \rg)=\ZZ_{q+1}\times\ZZ_{q+1}\rtimes S$.

\vskip 3mm
(5) $r\di (q-1)$, $r\ne 2$:
\vskip 3mm

 Followed from Table \ref{character}, there is only one conjugate classes subgroup of order $r$, where $r\mid q-1$ and $r\neq 2$. Let $\lg z_{q-1} \rg$ be a subgroup of $\ZZ_{\frac{q^2-1}{d}}\cong L$ $(\leq T)$ that has prime order $r$. Note that $\N_T(L)\cong\ZZ_{\frac{q^2-1}{d}}\rtimes\ZZ_2$. Thus, we obtain $\N_{T}(\lg z_{q-1} \rg)=\ZZ_{\frac{q^2-1}{d}}\rtimes\ZZ_2$ from the structure of the subgroups of $T$.
By Table \ref{character}, we have that $|\C_T(\lg z_{q-1}\rg)|=\frac{q^2-1}{d}$.
It follows from Proposition \ref{fieldconjugate} that $|\C_R(\lg z_{q-1} \rg)|=q^2-1$. Therefore, we have that $\N_R(\lg z_{q-1} \rg)=\ZZ_{q^2-1}\rtimes\ZZ_2$.

\vskip 3mm
(6) $r\di q^2-q+1$:
\vskip 3mm
 Let $\lg z_{q^2-q+1} \rg$ be a character subgroup of $\ZZ_{\frac{q^2-q+1}{d}}\cong A$ ($\leq T$) with order $r$. Note that there is only one conjugate classes subgroup in $\PSU(3,q)$ of order $\frac{q^2-q+1}{d}$ and $\N_T(A)=\ZZ_{\frac{q^2-q+1}{d}}\rtimes\ZZ_3$. Thus, we obtain $\N_{T}(\lg z_{q^2-q+1} \rg)\cong\ZZ_{\frac{q^2-q+1}{d}}\rtimes\ZZ_3$, which is maximal in $T$.
By Table \ref{character}, we have that $|\C_T(\lg z_{q^2-q+1}\rg)|=\frac{q^2-q+1}{d}$.
It follows from Proposition \ref{fieldconjugate} that $|\C_R(\lg z_{q^2-q+1} \rg)|=q^2-q+1$. Therefore, we have that $\N_R(\lg z_{q^2-q+1} \rg)=\ZZ_{q^2-q+1}\rtimes\ZZ_3$.
\qed

\begin{sidewaystable}[thp]
\caption{Subgroups of prime order $r$ in $\PSU(3,q)$: conjugacy classes }\label{sy}
\centering
\begin{tabular}{|c|c|c|c|}
\hline $r$&$p$ &$p$&$2$\\

  \hline
 \tiny {representation  } &
 $\begin{array} {l} z_0={\tiny\left(
  \begin{array}{ccc}
    1 & 0 & \epsilon \\
    0& 1 & 0\\
    0& 0 & 1 \\
  \end{array}
\right) } \\ $\tiny{ J=]1,1,1[ }$  \end{array}$

&$\begin{array} {l} z_l'={\tiny\left(
         \begin{array}{ccc}
           1 & a^l & b \\
           0 & 1 & -a^{lq} \\
           0 & 0 & 1 \\
         \end{array}
       \right) }\\  $\tiny{J=]1,1,1[,}$  \end{array} $ &  $\begin{array}{l} z_2={\tiny\left(
  \begin{array}{ccc}
    1 & 0 & 0 \\
    0& -1 & 0\\
    0& 0 & 1 \\
  \end{array}
\right)} \\ ${\tiny J=[1,1,1]}$ \end{array}$\\
\hline

$|\C_T(\lg z\rg)|$&${q^3(q+1)}/d$&$q^2$&$q(q^2-1)(q+1)/d$\\
\hline
$|\N_T(\lg z\rg)|$&${q^3(q+1)(p-1)}/{d}$&$q^2(p-1)$&${q(q^2-1)(q+1)}/{d}$\\
\hline
$|\C_R(\lg z\rg)|$&$q^3(q+1)$&$q^2$&$q(q^2-1)(q+1)$\\
\hline
$|\N_R(\lg z\rg)|$&$q^3(q+1)(p-1)$&$q^2(p-1)$&$q(q^2-1)(q+1)$\\
\hline
 {\tiny Remark}& {\tiny$\epsilon^q+\epsilon=0$,$\epsilon\in\FF_{q^2}$; \text{and} $\epsilon=1$ if $p=2$} & {\tiny \makecell{$0\le l \le d-1$; \\$a^{lq+l}+b+b^q=0$,$b\neq0$;\\ $|a|=3$ and $b=-\frac{1}{2}$ if $d=3$;\\
  $z_0'^{R}=z_1'^{R}=z_2'^{R}$, if $d=3$}} & \tiny{$q $ odd}\\
  \hline
\hline
 $r$& $3$&$3$&$r\mid q+1$\\ \hline
{\tiny representation} & $\begin{array}{l} z_3^A={\tiny [1, \omega,  \omega^2] }\\ ${\tiny J=[1,1,1]}$ \end{array}$

 &  $\begin{array}{l} z_3^B ={\tiny  [\rho^{\frac{q+1}{9}}, \rho^{\frac{q+1}{9}},  \rho^{\frac{-2(q+1)}9} ] } \\    ${\tiny J=[1,1,1]}$ \end{array}$

   & $\begin{array}{l}  z_{q+1}^A(a)=\tiny{ [x, x^a, x^{b}] } \\ ${\tiny J=[1,1,1]}$ \end{array}$     \\

\hline
$|\C_T(\lg z\rg)|$&$(q+1)^2$&${q(q+1)^2(q-1)}/{3}$&${(q+1)^2}/{d}$\\
\hline
$|\N_T(\lg z\rg)|$&$s(q+1)^2$&${q(q+1)^2(q-1)}/{3}$&${s(q+1)^2}/{d}$\\
\hline
$|\C_R(\lg z\rg)|$&$3(q+1)^2$&$q(q+1)^2(q-1)$&$(q+1)^2$\\
\hline
$|\N_R(\lg z\rg)|$&$3s(q+1)^2$&$q(q+1)^2(q-1)$&$s(q+1)^2$\\
\hline
{\tiny Remark}& {\tiny\makecell{$\omega\in\FF_{q^2}$, s=1,2\\$|\omega|=3$, $3\mid q+1$}}
& {\tiny\makecell{$\rho\in\FF_{q^2}$;\\$|\rho|=q+1$, $9\mid q+1$}} & {\tiny\makecell{$x\in\FF_{q^2}$,$1+a+b\equiv 0\pmod r$;\\$|x|=r$, $r\mid q+1$, $s\mid 6$;\\
$2\le a\le \frac {r-1}2$,
}} \\
 \hline
 \hline
 $r$& $r\mid q+1$&$r\mid q-1$ &$r\mid q^2-q+1$\\
 \hline

 {\tiny representation} & $\begin{array}{l}   z_{q+1}^B :={\tiny
   [ x, x,  x^{-2}]}\\  $ {\tiny J=[1,1,1] }$ \end{array}      $ & $\begin{array}{l} z_{q-1}:={\tiny [ \th^{\frac{q-1}{r}}, 1, \th^{-\frac{q-1}{r}}]}\\ ${\tiny J=]1,1,1[}$ \end{array}$

  & $z_{q^2-q+1}$\\

\hline
$|\C_T(\lg z\rg)|$&${q(q+1)^2(q-1)}/{d}$&${q^2-1}/{d}$&${q^2-q+1}/{d}$\\
\hline
$|\N_T(\lg z\rg)|$&${q(q+1)^2(q-1)}/{d}$&${2(q^2-1)}/{d}$&${3(q^2-q+1)}/{d}$\\
\hline
$|\C_R(\lg z\rg)|$&$q(q+1)^2(q-1)$&$q^2-1$&$q^2-q+1$\\
\hline
$|\N_R(\lg z\rg)|$&$q(q+1)^2(q-1)$&$2(q^2-1)$&$3(q^2-q+1)$\\
\hline

\tiny{ Remark}&{\tiny \makecell{$x\in\FF_{q^2}^*$;\\$|x|=r$, $r\mid q+1$}}& {\tiny\makecell{$\FF_{q}^*=\lg \th\rg$}} & * \\
 \hline

\end{tabular}

\end{sidewaystable}

\begin{sidewaystable}[thp]
\caption{Elements of order greater than 2 in $\PSU(3,q)$: conjugacy classes}
\label{character}
\centering
\begin{tabular}{|c|c|c|c|c|c|c|}
  \hline
Label& number of classes & Representative $g$ & Gram matrix J &$|g|$ & $|C_T(g)|$&parameters \\

\hline
 $\mathcal{C}_2$& 1 & {\tiny $\left(
         \begin{array}{ccc}
           1 & 0 & \epsilon \\
           0 & 1 & 0 \\
           0 & 0 & 1 \\
         \end{array}
       \right)$}
 & {\tiny $]1,1,1[$} &$p$ & $\frac{q^3(q+1)}{d}$ & \makecell{$\epsilon^q+\epsilon=0, \epsilon\in\FF_{q^2}$; \\
and $\epsilon=1$ if $p=2$ }\\
  \hline
$\mathcal{C}^{(l)}$&d&{\tiny $\left(
         \begin{array}{ccc}
           1 & a^l & b\\
           0 & 1 & -a^{lq} \\
           0 & 0 & 1 \\
         \end{array}
       \right) $} &{\tiny $]1,1,1[$}  &$p$     &$q^2$&\makecell{$0\le l \le d-1$, \\$a^{lq+l}+b+b^q=0$,$b\neq0$;\\
      $|a|=3$ and $b=-\frac{1}{2}$ if $d=3$}\\
\hline
$\mathcal{C}_4^{(k)}$&  $\frac{q+1}{d}-1$
&{\tiny $[\rho^k, \rho^k, \rho^{-2k}]$}  &{\tiny $[1,1,1]$}  &$|g|\mid \frac{q+1}k$   &$\frac{q(q+1)^2(q-1)}{d}$&\makecell{$1\leq k\leq \frac{q+1}{d}-1$;\\$\rho\in\FF_{q^2}$, $|\rho|=q+1$}\\

\hline
$\mathcal{C}_6'$&$\frac{d-1}{2}$&{\tiny $[
           1,\omega,\omega^2]$} & {\tiny $[1,1,1]$}  & 3   &$(q+1)^2$&\makecell{$\omega\in\FF_{q^2}$;\\$|\omega|=3$}\\
\hline
$\mathcal{C}_6^{(k,l,m)}$&$\frac{\frac{q^2-q+1}{d}-1}{6}$& {\tiny $
           [\rho^k, \rho^l, \rho^m]$}  & {\tiny $[1,1,1]$}   &$|g|\di (q+1)$   &$\frac{(q+1)^2}{d}$&\makecell{$1\leq k<l\leq \frac{q+1}{d}$;\\$l<m\leq q+1$;\\$k+l+m\equiv 0(\mod  q+1)$;\\$\rho\in\FF_{q^2}, |\rho|=q+1$}\\
\hline
$\mathcal{C}_7^{(k)}$&$\frac{\frac{q^2-q+1}{d}-1}{2}-\frac{3-d}{2}$
&{\tiny $[ \sigma^k, \rho^k,\sigma^{-qk}]$}  &{\tiny $]1,1,1[$}  &$|g|\di (q^2-1)$                &$\frac{q^2-1}{d}$&\makecell{$\rho,\sigma\in\FF_{q^2}$;\\ $|\rho|=q+1,|\sigma|=q^2-1$;\\
$1\leq k\leq \frac{q^2-1}{d}$;\\$k \not\equiv 0(\mod q-1)$;\\$\mathcal{C}_7^{(k)}=\mathcal{C}_7^{(-qk)}$}\\
\hline
$\mathcal{C}_8^{(k)}$&$\frac{\frac{q^2-q+1}{d}-1}{3}$& $*$ &$*$    &$|g|\mid q^2-q+1$ &$\frac{q^2-q+1}{d}$&\makecell{$\tau\in\FF_{q^6}$,$|\tau|=q^2-q+1$;\\
$1\leq k\leq \frac{q^2-q+1}{d}-1$;\\
$\mathcal{C}_8^{(k)}=\mathcal{C}_8^{(-qk)}=\mathcal{C}_8^{(q^2k)}$}\\
\hline
\end{tabular}
\end{sidewaystable}

\vskip 3mm

For $q' = 2$, the group $\PGU(3,2) \rtimes \langle f \rangle$ is soluble. In this case, the BG-Conjecture is established by \cite{BH}.
For $e=3$ and $q'=3$, one may confirm the BG-Conjecture by  Magma.  From now on, we assume that either $e\geq 5$ and $q'\geq 3$ or $e=3$ and $q'\geq 4$.
We will address the cases where $(\frac{q+1}{q'+1},3)$ equals $1$ and $3$ in Sections 5 and 6, respectively.

\section{$M_0\cong\PSU(3,q').(\frac{q+1}{q'+1},3)$ with $(\frac{q+1}{q'+1},3)=1$}\label{c=1}
Let  $q=p^{m}$, $q'=p^{m'}$ for some odd prime $e$ such that  $m=m'e$. Assume that $c:=(\frac{q+1}{q'+1},3)=1$.  Moreover, let $e\geq 5$ and $q'\geq 3$ or $e=3$ and $q'\geq 4$.
From the proof of Theorem~\ref{main1}, we get $(3,q+1)=:d=d':=(3,q'+1)$.
Let $${\small \begin{array}{ll} &T:=\PSU(3,q)\le G_1:=\PGU(3,q)=T\lg \d\rg \le G:=\PGammaU(3,q)= G_1\lg f\rg ;\\
&M_0:=\PSU(3, q')\le M_1:=\PGU(3,q')=M_0\lg \delta\rg \le M:=\PGU(3,q'):\lg f\rg=M_1 \lg f\rg; \\
&\O=[G:M]. \end{array} }$$
In what follows, we shall consider the action of $G$ on $\O=[G:M]$ and prove $\check{Q}(G)< \frac12$ so that the BG-Conjecture is true for $G$. By Remark~\ref{big}, this implies that the BG-Conjecture
is true for every group between $T$ and $G$ with the corresponding point-stabilizers.

 First we determine the representatives $\lg z\rg$ of  conjugacy classes of subgroups of prime order in $M$, while
 these subgroups are devided into two subclasses: those are contained in $M_1$ or those are not. We shall deal with them in Lemmas \ref{QM} and \ref{QF}, separately; measure $Q(\lg z \rg)$ for these subgroups $\lg z \rg$, where
\vspace{-5pt}$$ {\small Q(\lg z \rg):=\frac{|M|}{|\Omega|}\frac{(|z|-1)\Fix(\lg z \rg)}{|\N_{M}(\lg z\rg)|};}$$
 and then prove  $\check{Q}(G)< \frac12$ by Proposition~\ref{prob} so that  Theorem~\ref{saxl3} holds  for the case    $c=1$.
\begin{lem}\label{QM}
Let $c=1$ and  either $e\geq 5$ and $q'\geq 3$ or $e=3$ and $q'\geq 4$, where $q'=p^{m'}$.  Then using the notations given in Table 1 and Lemma \ref{PGU}, every cyclic subgroup  $\lg z \rg$ of prime order $r$ of $M_1$ is $M$-conjugate to  one of the following groups.

\vskip 3mm
{\rm(I)} Cases either $e\geq 5$ and $q'\geq 3$ or $e=3$ and $q'\geq 5:$

{\small\begin{enumerate}
\item[\rm(1)]$r=2$ and  $p$ is odd: $\lg z_2 \rg$ with
$Q(\lg z_2 \rg)<\frac{1}{26}.$
\item[\rm(2)]  $r=p$: $\lg z_0\rg$ with  $Q(\lg z_0 \rg)<\frac{1}{26}$; and  $\lg z_0'\rg$ with $Q(\lg z_0'\rg)<\frac{1}{26}.$
\item[\rm(3)] $r=3$ and  $r\mid q'+1$: $\lg z_A \rg$ with $Q(\lg z_A \rg)<\frac{1}{26};$
      $\lg z_B \rg$ with $Q(\lg z_B \rg)<\frac{1}{26};$ and $\lg \d'\rg\in M_1\setminus M_0$ with   $Q(\lg \delta'\rg)< \frac{1}{26}.$
\item[\rm(4)]  $r\mid q'-1$ and $r\neq 2$:  $\lg z_{q'-1} \rg$ with
$Q(\lg z_{q'-1} \rg)<\frac{1}{26(q'-1)}.$
\item[\rm(5)]  $r\mid q'+1$ and  $r\neq 2,3$: $\lg z_{q'+1}^A (a)\rg$ with  $Q(\lg z_{q'+1}^A(a) \rg)<\frac{1}{13(q'+1)^2};$ and $\lg z_{q'+1}^B \rg$ with
    $Q(\lg z_{q'+1}^B \rg)<\frac{1}{26(q'+1)}.$
\item[\rm(6)] $r\di q'^2-q'+1$ with $r\neq 3$: $\lg z\rg=\lg z_{q'^2-q'+1} \rg.$
    If $(e,p,m')\notin\{(3,2,2),(3,2,4),(3,3,2)\},$
    then $Q(\lg z_{q'^2-q'+1} \rg)<\frac{1}{26(q'^2-q'+1)}.$
    Otherwise, there is a unique conjugacy class of subgroups of order $r$, and we have $Q(\lg z_{q'^2-q'+1} \rg)<\frac{1}{26}.$
\end{enumerate}}
\vskip 3mm
{\rm(II)} Case $e=3$ and $q'=4:$

{\small\begin{enumerate}
\item[\rm(1)]  $r=2$: $\lg z_0\rg$ with  $Q(\lg z_0 \rg)<0.1\times\frac{1}{26}$; and  $\lg z_0'\rg$ with $Q(\lg z_0'\rg)<0.85\times\frac{1}{26}.$
\item[\rm(2)]  $r=3$:  $\lg z_{q'-1} \rg$ with
$Q(\lg z_{q'-1} \rg)<0.02\times\frac{1}{26}.$
\item[\rm(3)]  $r=5$: $\lg z_{q'+1}^A (a)\rg$ with  $Q(\lg z_{q'+1}^A(a) \rg)<0.75\times\frac{1}{26};$ and $\lg z_{q'+1}^B \rg$ with
    $Q(\lg z_{q'+1}^B \rg)<0.28\times\frac{1}{26}.$
\item[\rm(4)] $r=13$: $\lg z_{q'^2-q'+1} \rg$ and
$Q(\lg z_{q'^2-q'+1} \rg)<0.5\times\frac{1}{26}.$
\end{enumerate}}
\end{lem}
\demo  Suppose  $c=1$ and either $e\geq 5$ and $q'\geq 3$ or $e=3$ and $q'\geq 5$.
The  conclusion for  $e=3$ and $q'=4$, that is (II) in the lemma,  can be computed  directly by inserting the exact value of $e,q'$ in the proof of (I).
 So   we assume either $e\ge 5$ and $q'\ge 3$ or $e=3$ and $q'\ge 7$  (if $e=3$, then  $c=1$ gives $q'\ne 5$).
First we have
\vspace{-5pt}{\small$$\begin{array}{lll}
\frac{|M|}{|\Omega|}&=&\frac{|M|^2}{|G|}=\frac{2m|\PGU(3,q')|^2}{|\PGU(3,q)|}=\frac{2mq'^6(q'^3+1)^2(q'^2-1)^2}{q^3(q^3+1)(q^2-1)}\\
&\leq&\frac{2mq'^{16}}{q^8} \frac{ (1+q'^{-3})^2}{(1-q^{-2})} \frac{(1-q'^{-2})^2}{1+q^{-3}}\le  \frac{2mq'^{16}}{q^8} \cdot \frac{13}{12} \cdot 1
\leq \frac{13m q'^{16}}{6q^8}.
\end{array}$$}\vspace{-5pt}
 By Lemmas ~\ref{PGU} and~\ref{primeproperties}, we need to discuss the following cases, separately.
\vskip 3mm
(1)  Case $r=2$ and $p$ is odd: By Lemma~\ref{primeproperties}, we may set $z=z_2$, while  $\lg z_2 \rg^f=\lg z_2 \rg$.
It follows from Proposition \ref{man} and Table \ref{sy} that  $\Fix(\lg z_2\rg)=\frac{|\N_{G}(\lg z_2\rg)|}{|\N_{M}(\lg z_2 \rg)|}=\frac{q(q^2-1)(q+1)\times 2m}{q'(q'^2-1)(q'+1)\times 2m}=\frac{q(q^2-1)(q+1)}{q'(q'^2-1)(q'+1)}$, so that
\vspace{-5pt}{\small$$\begin{array}{lll}
Q(\lg z_2 \rg)&=&\frac{|M|}{|\Omega|}\times\frac{(r-1)\times\Fix(\lg z_2 \rg)}{|\N_{M}(\lg z_2 \rg)|}=\frac{|M|}{|\Omega|}\frac{q(q^2-1)(q+1)}{2m[q'(q'^2-1)(q'+1)]^2}\\
&\leq&\frac{13mq'^{16}}{6q^8}\frac{q^4(1+q^{-1})}{2mq'^8(1-q'^{-2})^2}\le  \frac{793}{576q'^{4(e-2)}}
\lvertneqq \frac{1}{26}.
\end{array}$$}\vspace{-5pt}
\quad(2) Case  $r=p$:  By Lemma~\ref{primeproperties}, we may set either $z=z_0$ or $z=z_l'$ where $0\le l\le d'-1$.

\vskip 3mm
Suppose that  $z=z_0$. Then   $|\N_{G_1}(\lg z_0 \rg)|=q^3(q+1)(p-1)$  and $|\N_{M_1}(\lg z_0\rg)|=q'^3(q'+1)(p-1)$. Now, we have
\vspace{-5pt}{\small$$\begin{array}{lll}
Q(\lg z_0 \rg)&=&\frac{|M|}{|\Omega|}\frac{(p-1)\Fix(\lg z_0 \rg)}{|\N_{M}(\lg z_0\rg)|}=\frac{|M|}{|\Omega|}\frac{(p-1)||\N_{G}(\lg z_0 \rg)|}{|\N_{M}(\lg z_0\rg)|^2}\leq
\frac{|M|}{|\Omega|}\frac{2m(p-1) q^3(q+1)(p-1)} {(q'^3(q'+1)(p-1))^2}
\\&\leq&\frac{13m^2q'^{16}}{3q^8}\frac{q^4(1+q^{-1})}{q'^8}
\leq \frac{3172m^2}{729q}\frac{1}{ q'^{3e-8}}\lvertneqq \frac{1}{26}.
\end{array}$$}
For the case $e=3$ and $q'=4$, $|\N_G(\lg z_0\rg)|=2m|\N_{G_1}(\lg z_0\rg)|$ and $|\N_M(\lg z_0\rg)|=2m|\N_{M_1}(\lg z_0\rg)|$, as $p=2$, where $z_0$ is of form in Table 1. So we get
 $$Q(\lg z_0\rg)\leq\frac{|M|}{|\O|}\frac{(p-1) q^3(q+1)(p-1)} {(q'^3(q'+1)(p-1))^2\times 2m}<0.1\times\frac{1}{26}.$$

Suppose that $z=z_l'$ where $0\le l\le d'-1$. If $d'=3$, then three elements $z_l'$ are  conjugate to each other in $M_1$ and so set  $z=z_0'$.
 This implies that
\vspace{-5pt}{\small$$\begin{array}{lll}
Q(\lg z_0' \rg)&=&\frac{|M|}{|\Omega|}\times\frac{(p-1)\Fix(\lg z_0' \rg)}{|\N_{M}(\lg z_0' \rg)|}=\frac{|M|}{|\Omega|}\frac{(p-1)|\N_{G}(\lg z_0' \rg)|}{|\N_{M}(\lg z_0' \rg)|^2}
\le \frac{13mq'^{16}}{6q^8} \frac{ 2m(p-1)q^2(p-1)}{(q'^2(p-1))^2}
\leq\frac{13m^2}{3q} \frac {1}{q'^{5e-12}}\\
&\leq&\frac{13}{81} \frac {1}{q'^{5e-12}}\lvertneqq \frac{1}{26}.
\end{array}$$}
\quad(3)  Case $r=3$ and  $r\mid q'+1$:  By Lemmas~\ref{PGU}, we take $z=z_A$, $z_B$ and $\d'$. Hence, $\lg z\rg^f=\lg z\rg.$
\vskip 3mm
$z=z_A$:   $|\N_{G}(\lg z_A\rg)|=6ms(q+1)^2$ and $|\N_{M}(\lg z_A\rg)|=6ms'(q'+1)^2$, where $s,s'\in\{1,2\}$.
Hence,
\vspace{-5pt}{\small$$\begin{array}{lll}
Q(\lg z_A \rg)&=&\frac{|M|}{|\Omega|}\frac{2\Fix(\lg z_A \rg)}{|\N_{M}(\lg z_0\rg)|}=\frac{|M|}{|\Omega|}\frac{2|\N_{G}(\lg z_A \rg)|}{|\N_{M}(\lg z_A\rg)|^2}\leq
 \frac{13mq'^{16}}{6q^8}\frac{12ms(q+1)^2} {(6ms'(q'+1)^2)^2}\le \frac{3172}{2187}\frac1 {q'^{6(e-2)}}\lvertneqq \frac 1{26}.
\end{array}$$}\vspace{-5pt}

$z=z_B$:  $|\N_{G}(\lg z_B\rg)|=2mq(q+1)^2(q-1)$ and $|\N_{M}(\lg z_B\rg)|=2mq'(q'+1)^2(q'-1)$. Therefore, we have
{\small$$\begin{array}{lll}
Q(\lg z_B \rg)&=&\frac{|M|}{|\Omega|}\times\frac{2\Fix(\lg z_B\rg)}{|\N_{M}(\lg z_B\rg)|}=\frac{|M|}{|\Omega|}\frac{2|\N_{G}(\lg z_B \rg)|}{|\N_{M}(\lg z_3^B\rg)|^2}
\le \frac{13mq'^{16}}{6q^8} \frac{4mq(q+1)^2(q-1)}{(2mq'(q'+1)^2(q'-1))^2} <  \frac 5{q'^{4(e-2)}}\lvertneqq \frac{1}{26}.
\end{array}$$}

\vspace{-5pt}
$z=\d'$:  Since $c=1$, we get that $\delta'=\delta$. Then $|\N_{G}(\lg \d' \rg)|=|(\ZZ_{q^2-q+1}:\ZZ_3)||\lg f\rg|=6m(q^2-q+1)$ and $|\N_{M}(\lg \d' \rg)|=6m(q'^2-q'+1)$. Consequently, we have
\vspace{-5pt}{\small$$\begin{array}{lll}
Q(\lg \d' \rg)&=&\frac{|M|}{|\Omega|}\times\frac{2\Fix(\lg \d' \rg)}{|\N_{M}(\lg \d' \rg)|}=\frac{|M|}{|\Omega|}\frac{2|\N_{G}(\lg \d' \rg)|}{|\N_{M}(\lg \d' \rg)|^2}
\le \frac{13mq'^{16}}{6q^8} \frac{12m(q^2-q+1)}{(6m(q'^2-q'+1))^2}
\le \frac{13}{8} \frac 1{q'^{6(e-2)}}\lvertneqq \frac{1}{26}.
\end{array}$$}
\vspace{-5pt}\quad(4) Case $r\mid q'-1$ and $r\neq 2$: By Lemma \ref{primeproperties}, we may set $z=z_{q'-1}$. Hence,
\vspace{-5pt}{\small$$
\begin{array}{lll}
(q'-1)Q(\lg z_{q'-1} \rg)&=&(q'-1)\times\frac{|M|}{|\Omega|}\times\frac{(r-1)\Fix(\lg z_{q'-1} \rg)}{|\N_M(\lg z_{q'-1} \rg)|}=(q'-1)\times\frac{|M|}{|\Omega|}\times\frac{(r-1)|\N_G(\lg z_{q'-1} \rg)|}{|\N_M(\lg z_{q'-1} \rg)|^2}\\
&=&\frac{(q'-1)|M|}{|\Omega|}\times\frac{(r-1)2(q^2-1)\times 2m}{[2(q'^2-1)\times 2m]^2}
\leq\frac{13mq'^{18}}{6q^8}\frac{q^2}{4mq'^2(1-q'^{-2})^2}
\leq\frac{351}{512 q'^{6e-16}}\lneqq\frac{1}{26}.
\end{array}$$}
\quad(5) Case $r\mid q'+1$ and $r\neq 2,3$:
By Lemma \ref{primeproperties}, we take $z=z_{q'+1}^A(a)$ and $z=z_{q'+1}^B$.
\vskip 3mm

$z=z_{q'+1}^A(a):$ By the proof in Lemma \ref{primeproperties}, we get that if  $z_{q'+1}^A(a_1)$ and $z_{q'+1}^A(a_2)$ are not conjugate in $M_1$ for some integers $a_1$ and $a_2$, then they are also not conjugate in $G_1$. Since $\lg z_{q'+1}^A(a) \rg^f=\lg z_{q'+1}^A(a) \rg$, by Lemma \ref{primeproperties}, we have $|\N_G(\lg z_{q'+1}^A(a)\rg)|=s(q+1)^2\times 2m$ and $|\N_M(\lg z_{q'+1}^A(a)\rg)|=s'(q'+1)^2\times 2m$, where $s,s'\mid 6$.
Therefore, we have
\vspace{-5pt}{\small$$\begin{array}{lll}
\frac{(q'+1)^2}{2}Q(\lg z_{q'+1}^A(a) \rg)&=&\frac{(q'+1)^2}{2}\times\frac{|M|}{|\Omega|}\times\frac{(r-1)\Fix(\lg z_{q'+1}^A(a) \rg)}{|\N_{M}(\lg z_{q'+1}^A(a)\rg)|}
=\frac{(q'+1)^2}{2}\times\frac{|M|}{|\Omega|}\frac{(r-1)s(q+1)^2}{s'(q'+1)^2s'(q'+1)^2\times 2m}\\
&\leq&\frac{13mq'^{16}}{12q^8} \frac{3(r-1)(q+1)^2}{m(q'+1)^2}
\leq\frac{23q'^{15}}{7q^6}\lneqq\frac{1}{26} ~~\text{(as $r-1\leq q'$)}.
\end{array}$$}
$\quad z=z_{q'+1}^B$: The number of conjugate classes of type $z_{q'+1}^B$ is only one. Since $\lg z_{q'+1}^B\rg^f=\lg z_{q'+1}^B\rg$, by Lemma \ref{primeproperties}, we have $|\N_G(\lg z_{q'+1}^B\rg)|=q(q-1)(q+1)^2\times 2m$ and $|\N_M(\lg z_{q'+1}^B\rg)|=q'(q'-1)(q'+1)^2\times 2m$.
Noting  $q'\geq 9$ and $r-1\leq q'$ in this case, we have
\vspace{-5pt}{\small$$
\begin{array}{lll}
(q'+1)Q(\lg z_{q'+1}^B \rg)&=&\frac{(q'+1)|M|}{|\Omega|}\times\frac{(r-1)|\Fix(z_{q'+1}^B)|}{\N_{M}(\lg z_{q'+1}^B\rg)}
=\frac{(q'+1)|M|}{|\Omega|}\times\frac{(r-1)q(q-1)(q+1)^2}{2m[q'(q'-1)(q'+1)^2]^2}\\
&\leq&\frac{13mq'^{16}}{6q^8}\times\frac{(r-1)q^4(1+q^{-1})^2}{2mq'^7(1-q'^{-1})^2}
\leq\frac{5q'^{10}}{2q^4}\lneqq\frac{1}{26}.
\end{array}$$}
\quad(6) Case $r\mid q'^2-q'+1$ and $r\neq 3:$
By Lemma \ref{primeproperties}, we may set $z=z_{q'^2-q'+1}$. Note that there is only one conjugate class subgroup with order $r$, and {\small$\Fix(\lg z_{q'^2-q'+1} \rg)=\frac{|\N_{G}(\lg z_{q'^2-q'+1}\rg)|}{|\N_{M}(\lg z_{q'^2-q'+1}\rg)|}$}, {\small$|\N_{G_1}(\lg z_{q'^2-q'+1}\rg)|=3(q^2-q+1)$} and {\small$|\N_{M_1}(\lg z_{q'^2-q'+1}\rg)|=3(q'^2-q'+1)$} from Proposition \ref{man} and Table 1. We now discuss the cases of $e=3$ and $e\geq 5$ separately.
\vskip 3mm
$e=3:$ Since $q=q'^3$ and $q'^2-q'+1\mid q^2-1$, the element $z_{q'^2-q+1}$ can be represented by a diagonal matrix of order $r$. Hence, we conclude that {\small$\Fix(\lg z_{q'^2-q'+1} \rg)=\frac{2m|\N_{G_1}(\lg z_{q'^2-q'+1}\rg)|}{2m|\N_{M_1}(\lg z_{q'^2-q'+1}\rg)|}=\frac{q^2-q+1}{q'^2-q'+1}$,} due to $\lg z_{q'^2-q'+1}\rg^f=\lg z_{q'^2-q'+1}\rg.$

If $e=3$ and $p=2$, then $m'\neq 3$, as $c=1$.
For $(e,p,m')\in\{(3,2,2),(3,2,4),(3,3,2)\}$, there is only one conjugacy classes subgroup of order $r\mid q'^2-q'+1$, as $q'^2-q'+1$ is a prime. This implies that
${\small\begin{array}{l}
Q(\lg z_{q'^2-q'+1} \rg)=\frac{|M|}{|\Omega|}\frac{(r-1)\Fix(\lg z_{q'^2-q'+1}\rg)}{|\N_{M}(\lg z_{q'^2-q'+1}\rg)|}\leq\frac{|M|}{|\Omega|}\frac{(r-1)(q^2-q+1)}{3(q'^2-q'+1)^2}\lneqq\frac{1}{26},
\end{array}}$
by inserting the these values of $e,p,m'$.
And for other cases, we have
\vspace{-5pt}{\small$$
\begin{array}{lll}
&&(q'^2-q'+1)Q(\lg z_{q'^2-q'+1} \rg)=\frac{(q'^2-q'+1)|M|}{|\Omega|}\times\frac{(r-1)\Fix(\lg z_{q'^2-q'+1}\rg)}{|\N_{M}(\lg z_{q'^2-q'+1}\rg)|}\\
&&\leq\frac{|M|(q'^2-q'+1)}{|\Omega|}\frac{(r-1)(q^2-q+1)}{3(q'^2-q'+1)^2}
\leq\frac{13mq'^{16}}{6q^8}\frac{(q^2-q+1)}{3}
\leq\frac{13m'e}{18p^{6m'e-16m'}}\lneqq\frac{1}{26}.
\end{array}$$}

$e\geq 5:$ We obtain that
\vspace{-5pt}{\small$$
\begin{array}{lll}
&&(q'^2-q'+1)Q(\lg z_{q'^2-q'+1} \rg)=\frac{(q'^2-q'+1)|M|}{|\Omega|}\times\frac{(r-1)\Fix(\lg z_{q'^2-q'+1}\rg)}{|\N_{M}(\lg z_{q'^2-q'+1}\rg)|}\\&&\leq\frac{|M|(q'^2-q'+1)}{|\Omega|}\frac{(r-1)(q^2-q+1)\times 2m}{3(q'^2-q'+1)^2}
\leq\frac{13mq'^{16}}{6q^8}\frac{(q^2-q+1)\times 2m}{3}
\leq\frac{26m'^2e^2}{18p^{6m'e-16m'}}\lneqq\frac{1}{26}.\hskip 2.5cm \Box
\end{array}$$}

\begin{lem}\label{QF} Let $c=1$ and  either $e\geq 5$ and $q'\geq 3$ or $e=3$ and $q'\geq 4$. Then every cyclic subgroup  $\lg z \rg$ of prime order $r$ of $M$  not contained in $M_1$
is $M$-conjugate to  one of the following groups.
\vskip 3mm
{\rm(I)} Cases either $e\geq 5$ and $q'\geq 3$ or $e=3$ and $q'\geq 5:$
{\small\begin{enumerate}
 \item[\rm(1)]  $r\notin \{e,2\}$: $z=f^{\frac{2m}r}$ with $Q(\lg f^{\frac{2m}r}\rg)< \frac{1}{26m'}.$
\item[\rm(2)] $r=2:$  $z=f^m$ with $Q(\lg f^m \rg)< \frac{1}{26}.$
\item[\rm(3)] $r=e:$
\begin{enumerate}
\item[\rm(3.1)]  $r\nmid |M_1|:$  $z=f^{\frac{2m}{r}}$ with $Q(z) < \frac 1{26};$
\item[\rm(3.2)]    $r=p:$   $z=z_0 f^{\frac{2m}{r}}$ or  $z_0'f^{\frac{2m}{r}}$ with $Q(z_0 f^{\frac{2m}{r}})+Q(z_0' f^{\frac{2m}{r}})< \frac 1{26};$
 \item[\rm(3.3)] $r\mid q'+1$ and $r\neq 3:$ $z=f^{\frac{2m}r}$,   $(z_{q'+1}^A(a))^if^{\frac{2m}r}$ where $2\le a\le \frac{r-1}2$ and $1\leq i\leq r-1$, or  $(z_{q'+1}^B)^i f^{\frac{2m}r}$ with $1\leq i\leq r-1$,  and
    $Q(\lg f^{\frac{2m}r}\rg )+\sum\limits_{i=1}^{r-1}\sum\limits_{a=2}^{\frac{r-1}2}Q(\lg z_{q'+1}^A(a))^if^{\frac{2m}r}\rg )+\sum\limits_{i=1}^{r-1}Q(\lg z_{q'+1}^Bf^{\frac{2m}r}\rg )< \frac 1{26};$
 \item[\rm(3.4)]    $r\mid q'-1:$  $z_{q'-1}^i f^{\frac{2m}r}$ where $0\le i\le r-1$, and $\sum\limits_{i=0}^{r-1} Q(\lg z_{q'-1})^i f^{\frac{2m}r}\rg )< \frac 1{26};$
 \item[\rm(3.5)] $r\di (q'^2-q'+1):$ $z=z_{q'^2-q'+1}^if^{\frac{2m}{r}}$, $0\le i\le r-1$ and
   $\sum\limits_{i=0}^{r-1} Q(\lg z_{q'^2-q'+1}^if^{\frac{2m}{r}}\rg)< \frac{1}{26}.$
\end{enumerate}
\end{enumerate}}

{\rm(II)} Case $e=3$ and $q'=4:$
{\small\begin{enumerate}
\item[\rm(1)] $r=2:$  $z=f^m$ with $Q(\lg f^m \rg)< \frac{1}{20}\times\frac{1}{26}.$
\item[\rm(2)] $r=3:$ $z_{q'-1}^i f^{\frac{2m}r}$ where $0\le i\le 2$, and $\sum\limits_{i=0}^{r-1} Q(\lg z_{q'-1})^i f^{\frac{2m}r}\rg )< 4\times\frac 1{26}.$
 \end{enumerate}}
\end{lem}
\demo  Let $u\in M\setminus M_1$. Then $z=u_1f_1$ for some element $u_1\in M_1$ and $f_1\in\lg f \rg$. Since $|z|=r$, we get   $|f_1|=r$ and so we may set $z=u_1f'$, where $f':=f^{\frac{2m}{r}}$. In what follows, we discuss three cases, separately.  Similar to last lemma, we assume either $e\ge 5$ and $q'\ge 3$ or $e=3$ and $q'\ge 7$, while  the  conclusion for  $e=3$ and $q'=4$, that is (II) in the lemma,  can be read from the proof of (I) directly.

\vskip 3mm
(1)  Case $r\notin\{ e,2\}:$
\vskip 3mm
Suppose $r\notin\{ e,2\}$.   Then by Proposition \ref{fieldconjugate}, we mays set $u_1=1$ so that  $z=f^{\frac{2m}r}$.
Set $m=m'e=m_1re$ so that $\Aut(\PSU(3,q'))=\PGU(3,q')\rtimes\lg f^e \rg$. By Proposition \ref{fieldconjugate}, there is only one $\PGU(3,q')$-class $\lg f^{2m_1e} \rg$, which is split into $(3,\frac{q'+1}{q'^{\frac{1}{r}}+1})$ distinct $\PSU(3,q')$-classes by \cite[Proposition $3.3.12$]{BGbook}.
 Since $\N_{M}(\lg f^{\frac{2m}{r}} \rg)=\PGU(3,q'^{\frac{1}{r}})\rtimes\lg f \rg$ and $\N_{G}(\lg f^{\frac{2m}{r}} \rg)=\PGU(3,q^{\frac{1}{r}})\rtimes\lg f \rg$, we have
\vspace{-5pt}{\small$$\begin{array}{lll}
&&{\small \Fix(\lg f^{\frac{2m}{r}} \rg)}=\frac{|\PGU(3,q^{\frac{1}{r}})\rtimes\lg f \rg|}{|\PGU(3,q'^{\frac{1}{r}})\rtimes\lg f \rg|}=\frac{p^{3m_1e}(p^{3m_1e}+1)(p^{2m_1e}-1)}{p^{3m_1}(p^{3m_1}+1)(p^{2m_1}-1)},\\
&&m'Q(\lg f^{\frac{2m}{r}} \rg)=\frac{m'|M|}{|\Omega|}\frac{(r-1)\Fix(\lg f^{\frac{2m}{r}} \rg)}{|\N_{M}(\lg f^{\frac{2m}{r}} \rg)|}
=\frac{m'|M|}{|\Omega|}\frac{(r-1)q^{\frac{3}{r}}(q^{\frac{3}{r}}+1)(q^{\frac{2}{r}}-1)}{[q'^{\frac{3}{r}}(q'^{\frac{3}{r}}+1)(q'^{\frac{2}{r}}-1)]^2\times2m}\\
&&\leq \frac{13m'mq'^{16}}{6q^8} \frac{2(r-1)q^{\frac{8}{r}}}{[q'^{\frac{8}{r}}(1-2^{-2})]^2\times2m} \leq\frac{104}{27q'^{\frac{16e-38}{3}}}\lneqq\frac{1}{26}\, \text{(as $r\geq 3$, $r-1\leq q'$, $m'\leq q'$ and $p\geq 2$)}.
\end{array}$$}

(2) Case $r=2:$
\vskip 3mm
By Proposition \ref{fieldconjugate} again,  we take $z=f^{m}$.  By \cite[Proposition $3.3.15$]{BGbook} we have $\C_{M_1}(\lg f^m \rg)\cong\PGO(3,q')$ and $\C_{G_1}(\lg f^m \rg)\cong\PGO(3,q)$.
Hence, we have
\vspace{-5pt}$${\small\begin{array}{lll}
{\small\Fix(\lg f^m \rg)}&=&\frac{|\N_{G}(\lg f^m \rg)|}{|\N_{M}(\lg f^m \rg)|}=\frac{|\PGO(3,q)\rtimes\lg f \rg|}{|\PGO(3,q')\rtimes\lg f \rg|}=\frac{q(q^2-1)}{q'(q'^2-1)},\\
Q(\lg f^m \rg)&=&\frac{|M|}{|\Omega|}\frac{\Fix(\lg f^m \rg)}{|\N_{M}(\lg f^m \rg)|}=\frac{|M|}{|\Omega|}\frac{q(q^2-1)}{2m[q'(q'^2-1)]^2}
\leq \frac{13mq'^{16}}{6q^8}\frac{q^3}{2mq'^6(1-q'^{-2})^2}
\leq\frac{351}{256q'^{5e-10}}\lneqq\frac{1}{26}.
\end{array}}$$
(3) Case $r=e:$
\vskip 3mm
Suppose $r=e$.  Then we may set $u=u_1f^{\frac{2m}e}=u_1f^{2m'}$, where $u_1\in M_1$. Now we have $[x, f^{2m'}]=1$ for any $x\in M_1$ and so $|u_1|=r$.
Recalling $|M_1|=q'^3(q'^3+1)(q'^2-1)$,  our prime  $r$ has  the following possibilities:  $r\nmid |M_1|$, $r=p$, $r\mid q'+1$, $r\mid q'-1$ and $r\mid q'^2-q'+1$, which are discussed separately.
\vskip 3mm
(3.1) $r\nmid |M_1|$
\vskip 3mm
Up to conjugate in $M$, set $z=f^{2m'}$.
Hence, $\N_{G}(\lg z \rg)=\N_{M}(\lg z \rg)=M$ and
\vspace{-5pt}{\small$$\begin{array}{lll}
&& Q(\lg z \rg)=\frac{|M|}{|\Omega|}\frac{(r-1)\times 1}{|M|}
\leq \frac{13mq'^{16}}{6q^8}\frac{e}{2mq'^3(q'^3+1)(q'^2-1)}
\leq\frac{39eq'^8}{32q^{8}}\lneqq\frac{1}{26} \, \text{(as  $r-1\leq e$)}.
\end{array}$$}
\hskip 5mm(3.2) $r=p$
\vskip 3mm
Since there are two  $M_1$-conjugacy classes of elements of order $p$, with the representatives  $ z_0$ and $ z_0'$, we set $\lg z\rg=\lg f^{2m'}\rg$,
 $\lg z_0 f^{2m'}\rg$ and $\lg z_0'f^{2m'}\rg$.
By Proposition \ref{fieldconjugate}, both of them are conjugate to $\lg f^{2m'} \rg$ in $G_1$.
This implies that $|\N_{G}(\lg z_0f^{2m'}\rg)|=|\N_{G}(\lg z_0'f^{2m'}\rg)|=|\N_{G}(\lg f^{2m'}\rg)|=|M|$.
By Lemma \ref{primeproperties}, we conclude that $|\C_{M_1}(\lg z_0\rg)|=q'^3(q'+1)$ and $|\C_{M_1}(\lg z_0'\rg)|$
$=q'^2$. It deduces that
$|\N_{M}(\lg z_0f^{2m'}\rg)|=|\C_{M}(\lg z_0\rg)|\geq eq'^3(q'+1)$ and $|\N_{M}(\lg z_0'f^{2m'}\rg)|=|\C_{M}(\lg z_0'\rg)|\geq eq'^2$.
Therefore, we have
\vspace{-5pt}{\small
$$\begin{array}{ll}
&\begin{array}{ll} &\Fix(\lg z_0 f^{2m'} \rg)=\Fix(\lg z_0'f^{2m'}\rg)
=\Fix(\lg f^{2m'}\rg)\\
&=\frac{|\N_{G}(\lg z_0f^{2m'} \rg)|}{|\N_{M}(\lg z_0f^{2m'} \rg)|}+\frac{|\N_{G}(\lg z_0'f^{2m'} \rg)|}{|\N_{M}(\lg z_0'f^{2m'}\rg)|}+\frac{|\N_{G}(\lg f^{2m'} \rg)|}{|\N_{M}(\lg f^{2m'}\rg)|}\\
&\leq \frac{|M|}{eq'^3(q'+1)}
+\frac{|M|}{eq'^2}+\frac{|M|}{|M|}
=2m'(q'^3+1)(q'-1)[1+q'(q'+1)]+1,\end{array}\\
&\begin{array}{ll}
&Q(\lg f^{2m'}\rg)=\frac{|M|}{|\Omega|}\frac{(p-1)\Fix(\lg f'\rg)}{|\N_{M}(\lg f'\rg)|}\leq\frac{|M|}{|\Omega|}
(\frac{1+q'^2+q'}{q'^4 }+\frac{p-1}{2mq'^6(q'^2-1)})
\leq\frac{|M|}{|\Omega|}(\frac{5}{3q'^2}+\frac{9}{16mq'^7}),\\
\end{array}\\
&\begin{array}{ll}
&Q(\lg z_0f^{2m'}\rg)=\frac{|M|}{|\Omega|}\frac{(p-1)\Fix(\lg z_0f'\rg)}{|\N_{M}(\lg z_0f'\rg)|}\leq\frac{|M|}{|\Omega|}(\frac{2m'(q'^2-q'+1)(1+q'^2+q')}{q'^2}+\frac{1}{q'^4})
\leq\frac{|M|}{|\Omega|}{\tiny(10m'q'^2+\frac{1}{q'^4})}\\
\end{array}\\
&\begin{array}{ll}
&Q(\lg z_0'f^{2m'}\rg)=\frac{|M|}{|\Omega|}\frac{(p-1)\Fix(\lg z_0'f'\rg)}{|\N_{M}(\lg z_0'f'\rg)|}\leq\frac{|M|}{|\Omega|}(\frac{2m'(q'^3+1)(d'+q'^2+q')}{q' }+\frac{1}{q'^2})\leq\frac{|M|}{|\Omega|}(\frac{280m'q'^4}{81}+\frac{1}{q'^2})\\
\end{array},\\
&\begin{array}{ll}
&Q(\lg f'\rg)+Q(\lg z_0f' \rg)+Q(\lg z_0'f' \rg)\leq\frac{|M|}{|\Omega|}(\frac{8}{3q'^2}+\frac{1}{q'^4}+\frac{9}{16mq'^7}+10m'q'^2+\frac{280m'q'^4}{81})\\
&\quad \leq\frac{52m'pq'^{14}}{9q^8}+\frac{13m'pq'^{12}}{2q^8}+\frac{39q'^{9}}{32q^8}+\frac{65m'^2pq'^{18}}{3q^8}+\frac{1820m'^2pq'^{20}}{243q^8}\\
&\quad\leq\frac{52q'^{15}}{9q^8}+\frac{13q'^{13}}{2q^8}+\frac{39q'^{9}}{32q^8}+\frac{65q'^{19}}{3q^8}+\frac{1820q'^{21}}{243q^8}
\lneqq\frac{1}{26} \,\,(\text{as}\,\,r=e=p\geq3).
\end{array} \end{array} $$
}
\hskip 5mm(3.3) $r\mid q'+1$
\vskip 3mm
If $e=3$, then $r=3$ and so $3\di (q'+1)$. Thus, $c=(\frac{q+1}{q'+1},3)=(q'^2-q'+1, 3)=3$, a contradiction. Therefore,  $e\geq 5$.

Let $z=u_0f^{\frac{2m}r}$, where $u_0\in M_1$.  Since  every subgroup of order $r$ where $r\di q'+1$ of $M_1$  is conjugate to $\lg z_{q'+1}^A(a)\rg$ where $2\le a\le \frac{r-1}2$ or  $\lg z_{q'+1}^B\rg$, we may set
  $\lg z\rg=\lg f^{\frac{2m}r}\rg$,  $\lg (z_{q'+1}^A(a))^if^{\frac{2m}r}\rg$  or $\lg (z_{q'+1}^B)^i f^{\frac{2m}r}\rg$, where $1\leq i\leq r-1$.
  By Proposition \ref{fieldconjugate}, all the elements $u_0f^{\frac{2m}r}$ of  order $r$ in $M$  are $G_1$-conjugate to $ f^{\frac{2m}r}$. Hence, $|\N_{G}(\lg (z_{q'+1}^A(a))^if'\rg)|=|\N_{G}(\lg (z_{q'+1}^B)^if'\rg)|=|\N_{G}(\lg f'\rg)|=|M|$. By Table \ref{sy}, we get that \vspace{-5pt}$${\small\begin{array}{ll}
  &|\N_{M}(\lg (z_{q'+1}^A(a))^if'\rg)|=|\C_{M}(\lg (z_{q'+1}^A(a))^i \rg)|=|\C_{M_1}(\lg (z_{q'+1}^A(a))^i\rg)|e=(q'+1)^2e,\\
&|\N_{M}(\lg (z_{q'+1}^B)^if'\rg)|=|\C_{M}(\lg (z_{q'+1}^B)^i\rg)|=|\C_{M_1}(\lg (z_{q'+1}^B)^i\rg)|e=q'(q'+1)^2(q'-1)e.
\end{array}}$$\vspace{-5pt}
Now, we conclude that
\vspace{-5pt}{\small$$ \begin{array}{ll}
&\begin{array}{l}
\Fix(\lg f'\rg)=\Fix(\lg (z_{q'+1}^A(a))^if' \rg)=\Fix(\lg (z_{q'+1}^B)^if'\rg)\\
\qquad \leq \frac{(r-3)(r-1)}{2}\frac{|\N_{G}(\lg (z_{q'+1}^A(a))^if'\rg)|}{|\N_{M}((\lg z_{q'+1}^A(a))^i f' \rg)|}+(r-1)\frac{|\N_{G}(\lg (z_{q'+1}^B)^if'\rg)|}{|\N_{M}(\lg (z_{q'+1}^B)^i f' \rg)|}
+\frac{|\N_{G}(\lg f'\rg)|}{|\N_{M}(\lg f'\rg)|}\\
\qquad\leq\frac{(r-3)(r-1)}{2}\frac{|M|}{(q'+1)^2e}
+(r-1)\frac{|M|}{q'(q'+1)^2(q'-1)e}
+\frac{|M|}{|M|}
 \leq emq'^6+2mq'^4+1;
\end{array}\\
&\begin{array}{l}
Q(\lg f' \rg)=\frac{|M|}{|\Omega|}\frac{(r-1)\Fix(\lg f' \rg)}{\N_{M}(\lg f'\rg)} \leq\frac{13mq'^{16}}{6q^8}\frac{r-1}{2mq'^3(q'^3+1)(q'^2-1)}
(emq'^6+2mq'^4+1)\\
\qquad \leq \frac{39e^2mq'^{14}}{32q^8}+\frac{39emq'^{12}}{16q^8}
+\frac{39eq'^8}{32q^8};
\end{array}\\
&\begin{array}{l}
 Q(\lg (z_{q'+1}^A(a))^if' \rg)=\frac{|M|}{|\Omega|}
 \frac{(r-1)\Fix(\lg (z_{q'+1}^A(a))^if'\rg)}{\N_{M}(\lg (z_{q'+1}^A(a))^if'\rg)}
\leq\frac{13mq'^{16}}{6q^8}\frac{r-1}{(q'+1)^2e}
(emq'^6+2mq'^4+1)\\
\qquad \leq \frac{13e m^2q'^{20}}{6q^8}+\frac{13m^2q'^{18}}{3q^8}
+\frac{13mq'^{14}}{6q^8};
\end{array}\\
&\begin{array}{l}
Q(\lg (z_{q'+1}^B)^if' \rg)=\frac{|M|}{|\Omega|}\frac{(r-1)\Fix(\lg (z_{q'+1}^B)^if' \rg)}{\N_{M}(\lg (z_{q'+1}^B)^if'\rg)}
\leq\frac{13mq'^{16}}{6q^8}\frac{r-1}{q'(q'+1)^2(q'-1)e}
(emq'^6+2mq'^4+1)\\
\qquad \leq
\frac{13e m^2q'^{18}}{4q^8}+\frac{13m^2q'^{16}}{2q^8}
+\frac{13mq'^{12}}{4q^8};
\end{array}\\
&\begin{array}{l}
Q(\lg f' \rg)+\sum\limits_{i=1}^{r-1}\sum\limits_{a=2}^{\frac{r-1}{2}}Q(\lg (z_{q'+1}^A(a))^if' \rg)+\sum\limits_{i=1}^{r-1}Q(\lg (z_{q'+1}^B)^if' \rg)
\leq \frac{39e^2mq'^{14}}{32q^8}+\frac{39emq'^{12}}{16q^8}
+\frac{39eq'^8}{32q^8}\\
\qquad +\frac{(r-1)(r-3)}{2}(\frac{13e m^2q'^{20}}{6q^8}+\frac{13m^2q'^{18}}{3q^8}
+\frac{13mq'^{14}}{6q^8})
+(r-1)(\frac{13e m^2q'^{18}}{4q^8}+\frac{13m^2q'^{16}}{2q^8}
+\frac{13mq'^{12}}{4q^8})\\\end{array}\end{array}$$}
{\small$$\begin{array}{l}
\qquad\leq
\frac{221e^2mq'^{14}}{96q^8}+\frac{91emq'^{12}}{16q^8}
+\frac{65e^2m^2q'^{18}}{12q^8}
+\frac{13e^3m^2q'^{20}}{12q^8}
+\frac{13e m^2q'^{16}}{2q^8}
+\frac{39eq'^{8}}{32q^8}\\
\qquad\leq
\frac{221\times 5^3q'^{14}}{96\times 2^5q^7}+\frac{91\times5^2q'^{12}}{16\times 2^5q^7}
+\frac{65\times 5^4q'^{18}}{12\times 2^5q^7}
+\frac{13\times 5^5q'^{20}}{12\times 2^5q^7}
+\frac{13\times5^3q'^{16}}{2^6q^7}
+\frac{39\times 5q'^{8}}{32\times 2^5q^7}\lneqq\frac{1}{26}.
\end{array}
$$}
\hskip 4mm(3.4) $r\mid q'-1$
\vskip 3mm

Let $z=u_0f^{\frac{2m}r}$, where $u_0\in M_1$.  Since  every subgroup of order $r$ where $r\di q'-1$ of $M_1$  is conjugate to  $\lg z_{q'-1}\rg$, we set $\lg z\rg=\lg z_{q'-1}^if^{\frac{2m}r}\rg$ where $0\le i\le r-1$.
By Proposition \ref{fieldconjugate},   all the elements $uf'\in M$ with order $r$ are $G_1$-conjugate to $f'$. Hence, $|\N_{G}(\lg z_{q'-1}^i f^{\frac{2m}{r}}\rg)|=|\N_{G}(\lg f^{\frac{2m}{r}}\rg)|=|M|$. By the arguments above and Table \ref{sy}, we get that $|\N_{M}(\lg z_{q'-1}^if'\rg)|=|\C_{M}(\lg z_{q'-1}\rg)|=2e|\C_{M_1}(\lg z_{q'-1}\rg)|=2e(q'^2-1)$, as $z_{q'-1}\in\{ h(k)\mid k\in\FF_{q'}\}$.
Hence,
\vspace{-5pt}{\small$$
\begin{array}{ll}
&\begin{array}{l}
\Fix(\lg f'\rg)=\Fix(\lg z_{q'-1}^if' \rg)\leq\frac{|\N_{G}(\lg f'\rg)|}{|\N_{M}(\lg f'\rg)|}+\sum\limits_{i=1}^{r-1}\frac{|\N_{G}(\lg f'\rg)|}{|\N_{M}(\lg z_{q'-1}^if'\rg)|}\\
\qquad\leq \frac{|M|}{|M|}+(r-1)\frac{|M|}{2e(q'^2-1)}
=1+(e-1)\frac{2mq'^3(q'^3+1)(q'^2-1)}{2e(q'^2-1)}\leq 1+mq'^3(q'^3+1);
\end{array}\\\end{array}$$}
{\small$$\begin{array}{ll}
&\begin{array}{l}
Q(\lg f' \rg)=\frac{|M|}{|\Omega|}\frac{(r-1)\Fix(\lg f' \rg)}{|\N_{M}(\lg f'\rg)|}
\leq
\frac{13mq'^{16}}{6q^8}\frac{(e-1)[1+mq'^3(q'^3+1)]}{2mq'^3(q'^3+1)(q'^2-1)}
\leq\frac{39eq'^8}{32q^8}+\frac{39meq'^{14}}{32q^8};
\end{array}\\
&\begin{array}{l}
Q(\lg z_{q'-1}^if' \rg)=\frac{|M|}{|\Omega|}\frac{(r-1)\Fix(\lg z_{q'-1}^if'\rg)}{|\N_{M}(\lg z_{q'-1}^if'\rg)|}
\leq\frac{13mq'^{16}}{6q^8}\frac{(r-1)[1+mq'^3(q'^3+1)]}{2e(q'^2-1)}
\leq\frac{39mq'^{14}}{32q^8}+\frac{91m^2q'^{20}}{72q^8};
\end{array}\\
&\begin{array}{l}
\sum\limits_{i=0}^{r-1}Q(\lg z_{q'-1}^if' \rg)
\leq \frac{39q'^9}{32q^8}+\frac{78emq'^{14}}{32q^8}
+\frac{91e m^2q'^{20}}{72q^8}
\lneqq\frac{1}{26}.
\end{array} \end{array} $$}
\hskip 4mm(3.5) $r\mid q'^2-q'+1$
\vskip 3mm
 If $e=3$, then it follows from $1=c=(\frac{q+1}{q'+1},3)=(q'^2-q'+1, 3)$ that $r=e\ne 3$. Therefore,  $e\geq 5$.

Set  $z=u_0f^{\frac{2m}{r}}$, where $u_0\in M_1$. Since  $M_1$ has only one class of subgroups $\ZZ_{q'^2-q'+1}$, we may set $u_0=z_{q'^2-q'+1}^i$, where $0\le i\le r-1$.
By Proposition \ref{fieldconjugate},   all these elements $z_{q'^2-q'+1}^if^{2m'}\in M$   are  $G_1$-conjugate to $f'$. Hence, $|\N_{G}(\lg z_{q'^2-q'+1}^if^{2m'}\rg)|=|\N_{G}(\lg f^{2m'}\rg)|=|M|$. By Table \ref{sy}, we get that

${\small|\N_{M}(\lg z_{q'^2-q'+1}^if^{2m'}\rg)|=|\C_{M}(\lg z_{q'^2-q'+1}^i\rg)|=|\C_{M_1}(\lg z_{q'^2-q'+1}^i\rg)|e}$
${\small=(q'^2-q'+1)e.}$

\f Therefore, we have
\vspace{-5pt}{\small$$\begin{array}{ll} &\begin{array}{ll}
&\Fix(\lg f^{2m'}\rg)=\Fix(\lg z_{q'^2-q'+1}^if^{2m'}\rg)\leq\frac{|\N_{G}(\lg f'\rg)|}{|\N_{M}(\lg f'\rg)|}+(r-1)\frac{|\N_{G}(\lg f^{2m'}\rg)|}{|\N_{M}(\lg z_{q'^2-q'+1}^if^{2m'}\rg)|}\\
&\qquad \leq \frac{|M|}{|M|}+(r-1)\frac{|M|}{(q'^2-q'+1)e}=1+(e-1)\frac{2mq'^3(q'^3+1)(q'^2-1)}{(q'^2-q'+1)e}\leq 1+2mq'^3(q'+1)(q'^2-1),
\end{array}\\
&\begin{array}{ll}
&Q(\lg f' \rg)=\frac{|M|}{|\Omega|}\frac{(r-1)|\Fix(\lg f' \rg)|}{|\N_{M}(\lg f'\rg)|}
\leq \frac{13mq'^{16}}{6q^8}\frac{(r-1)[1+2mq'^3(q'+1)(q'^2-1)]}{|M|}
\leq\frac{13mq'^{16}}{6q^8} (\frac{9}{16m'q'^8}+\frac{3e}{2q'^2})\\
&\qquad =\frac{39eq'^8}{32q^8}+\frac{13meq'^{14}}{4q^8},\\
\end{array}\\
&\begin{array}{ll}
&Q(\lg z_{q'^2-q'+1}^if' \rg)=\frac{|M|}{|\Omega|}\frac{(r-1)|\Fix(\lg z_{q'^2-q'+1}^if'\rg)|}{|\N_{M}(\lg z_{q'^2-q'+1}^if'\rg)|}
\leq\frac{13mq'^{16}}{6q^8}\frac{(r-1)[1+2mq'^3(q'+1)(q'^2-1)]}{(q'^2-q'+1)e}\\
 &\qquad=\frac{13mq'^{16}}{6q^8}(\frac{1}{q'^2(1-q'^{-1})}+\frac{2mq'^4(1+q'^{-1})^2}{q'^2(1-q'^{-1})})
\leq\frac{13mq'^{14}}{4q^8}+\frac{104m^2q'^{18}}{9q^8},
\end{array}\\
&\begin{array}{ll}
&
\sum\limits_{i=0}^{r-1}Q(\lg z_{q'^2-q'+1}^if' \rg)
\leq \frac{39eq'^8}{32q^8}+\frac{13meq'^{14}}{4q^8}
+(r-1)(\frac{13mq'^{14}}{4q^8}+\frac{104m^2q'^{18}}{9q^8})\\
&\qquad \leq
\frac{39eq'^8}{32q^8}+\frac{13emq'^{14}}{4q^8}
+\frac{13emq'^{14}}{4q^8}+\frac{104e m^2q'^{18}}{9q^8}\lneqq\frac{1}{26}.\hskip 6.2cm \Box
\end{array} \end{array}$$}

\begin{theorem}\label{saxl3}
With $c=(\frac{q+1}{q'+1},3)=1$, $e\geq 3$ an odd prime, $G=\PGammaU(3,q)$, and  $M=\PGU(3,q'):\lg f \rg$, the BG-Conjecture holds for the action of $G$ on $[G:M]$.
\end{theorem}
\demo Taking into account, all conjugate classes of subgroups $\lg z\rg $  of prime order in $M$ and an upper bound for the number $Q(\lg z\rg )$ have been determined  in  Lemma~\ref{QM} and Lemma \ref{QF},
  we are ready to show $\check{Q}(G)< \frac 12$, so that  the BG-Conjecture is true.
\vskip 3mm
Write $m'=r_1^{k_1} r_2^{k_2}\cdots r_{l}^{k_l}$ if $m'\neq 1$,  where $2\le r_1\lneqq r_2\lneqq \cdots \lneqq r_l$ are distinct primes and $k_i\geq 1$ for all $i$.
Further, let
\begin{enumerate}
  \item  [{\rm(a)}] $q'-1=2^{\epsilon_2}t_1^{j_1} t_2^{j_2}\cdots t_{o}^{j_o}$, where $3\le t_1\lneqq t_2\lneqq \cdots \lneqq t_o$ are distinct primes;
    \item  [{\rm(b)}] $q'+1=2^{\epsilon_1}3^{\kappa_1}s_1^{i_1} s_2^{i_2}\cdots s_{h}^{i_h}$,
 where $5\le s_1\lneqq s_1\lneqq \cdots \lneqq s_h$ are distinct primes; and
  \item [{\rm(c)}] $q'^2-q'+1=2^{\epsilon_3}3^{\kappa_3}u_1^{\ell_1} u_2^{\ell_2}\cdots u_{\imath}^{\ell_{\imath}}$ where $5\le \ell_1\lneqq \ell_2\lneqq \cdots \lneqq \ell_\imath$ are distinct primes.
\end{enumerate}
Furthermore, denote by $z_{q'-1, i}$ an element of order $t_i$ in $\lg z_{q'-1}\rg$; denote by $z_{q'+1, i}^A(a)$ and $z_{q'+1, i}^B$ elements of order $s_i$ in $\lg z_{q'+1}^A(a)\rg$ and $\lg z_{q'+1}^B\rg$, respectively; denote
 by $z_{q'^2-q'+1, i}$ an element of order $u_i$ in $\lg z_{q'^2-q'+1}\rg$.
Recall $\mathcal{P}^*(M):=\{\lg g_1\rg,\lg g_2\rg,\cdots,\lg g_k\rg\}$ is the set of representatives  of conjugacy classes of subgroups of prime order in $M$.

Suppose that either $e\ge 5$ and $q'\ge 3$ or $e=3$ and $q'\ge 7$. Then by Lemma \ref{QM}.(I) and Lemma \ref{QF}.(I), we get
$${\small\begin{array}{ll}
 &\check{Q}(G)=\frac{|M|}{|\O|}(\sum_{i=1}^k(|\lg g_i\rg|-1)\frac{\Fix(\lg g_i\rg)}{|\N_M(\lg g_i\rg )|})=\sum_{i=1}^kQ(\lg g_i\rg)\\
 &\qquad \leq (Q(\lg z_0 \rg)+Q(\lg z_0'\rg)
+Q(\lg z_2 \rg)+Q(\lg z_A\rg)
+Q(\lg z_B \rg)+Q(\lg \delta'\rg))\\
&\qquad+\sum\limits_{i=1}^{o}Q(\lg z_{q'-1, i}\rg)+
\sum\limits_{i=1}^{h}\sum\limits_{1\le a\le \frac{s_i-3}2} Q(\lg z_{q'+1, i}^A(a)\rg)+\sum\limits_{i=1}^{h}Q(\lg z_{q'+1, i}^B\rg)
+\sum\limits_{i=1}^{\imath} Q(\lg z_{q'^2-q'+1, i}\rg)\\
&\qquad+\sum\limits_{r_i\mid m,r_i\neq 2, e}Q(\lg f^{\frac{2m}{r_i}} \rg)
+Q(\lg f^m \rg)+\sum\limits_{  z\, {\rm  is\, in\, Lemma~\ref{QF}(I)(3)},  |z|=e}Q(\lg z\rg)\\
&\qquad\le  \frac{6}{26}+(q'-1) Q(\lg z_{q'-1}\rg)+ \frac{(q'+1)^2}{2} Q(\lg z_{q'+1}^A(a)\rg)+(q'+1)Q(\lg z_{q'+1}^B\rg)\\
&\qquad+(q'^2-q'+1)Q(\lg z_{q'^2-q'+1}\rg)+m'Q(\lg f^{\frac{2m}{r}}\rg)+\frac 1{26}+\frac 1{26}  \\
&\qquad<  \frac{6}{26}+7\cdot\frac{1}{26}=\frac 12.
\end{array}}$$
Suppose that  $e=3$ and $q'=4$  Then by Lemma \ref{QM}.(II) and Lemma \ref{QF}.(II), we have
\vspace{-5pt}$${\small\begin{array}{lll}
 \check{Q}(G) &\leq & (Q(\lg z_0 \rg)+Q(\lg z_0'\rg)
 +Q(\lg z_{q'-1, 1}\rg)+Q(\lg z_{q'+1, 1}^A(a)\rg)+Q(\lg z_{q'+1, 1}^B\rg)\\
&&+ Q(\lg z_{q'^2-q'+1, 1}\rg)
+Q(\lg f^m \rg)+\sum\limits_{  z\, {\rm  is\, in\, Lemma~\ref{QF}(II)(3)},  |z|=3}Q(\lg z\rg)< \frac 12,
\end{array}}$$
as desired. \qed
\section{$M_0\cong\PSU(3,q').(\frac{q+1}{q'+1},3)$ with $(\frac{q+1}{q'+1},3)=3$}\label{c=3}

Let  $q=p^{m}$, $q'=p^{m'}$ for some odd prime $e$ such that  $m=m'e$.  Assume that $c=(\frac{q+1}{q'+1},3)=3$. From the proof of Lemma~\ref{main1}, we get $d:=(3,q+1)=d':=(3,q'+1)=e=3$.
Let $${\small\begin{array}{ll} &T:=\PSU(3,q)\le G:=\PSigmaU(3,q)= T\lg f\rg ;\\
&M_0:=\PSU(3, q')\le M_1:=\PGU(3,q')=M_0\lg \delta\rg \le M:=\PGU(3,q'):\lg f\rg=M_1 \lg f\rg ;\\
&\O=[G:M]. \end{array} }$$
In what follows, we shall consider the action of $G$ on $\O$ and prove $\hat{Q}(G)< \frac12$, so that the BG-Conjecture is true for $G$. By Remark~\ref{big}, the BG-Conjecture
is true for  every  group $G_1$ between $T$ and $\PSigmaU(3,q)$ with the corresponding point-stabilizers.

 First we    determine the representatives $\lg z\rg$ of  conjugacy classes of subgroups of prime order in $M$, while
 these subgroups are devided into two subclasses: those are   contained $M_1$  or those not. We shall deal with them in Lemmas \ref{3SU} and \ref{3SUF}, separately; and also   measure $Q(\lg z \rg)$ for these subgroups $\lg z \rg$, where
\vspace{-5pt}$${\small Q(\lg z \rg):=\frac{|M|}{|\Omega|}\frac{(|z|-1)\Fix(\lg z \rg)}{|\N_{M}(\lg z\rg)|}};$$
and then prove the conjecture for this case in Theorem~\ref{c=3proof}.
Since $q' \neq 3, 4$ when $e = 3$, we may (and will) assume $q' \geq 5$ in the following analysis.

\begin{lem}\label{3SU}
Let $c=e=3$ and $q'\geq 5$. Then
using the notations given in Table \ref{sy} and Lemma \ref{PGU}, every cyclic subgroup $\lg z\rg$ of prime order $r$ of $M_1$ is $M$ conjugate to one of the following groups.
{\small\begin{enumerate}
\item[\rm(1)]$r=2$ and  $p$ is odd: $\lg z_2 \rg$ with
$Q(\lg z_2 \rg)<\frac{1}{26}.$
\item[\rm(2)]  $r=p$: $\lg z_0\rg$ with  $Q(\lg z_0 \rg)<\frac{1}{26}$; and  $\lg z_0'\rg$ with $Q(\lg z_0'\rg)<\frac{1}{26}.$
\item[\rm(3)] $r=3$ and  $r\mid q'+1$: $\lg z\rg=\lg z_A\rg$, $\lg z_B\rg$  or $\lg\d'\rg\in M_1\setminus M_0$ with $Q(\lg z_A\rg)+Q(\lg z_B\rg)+Q(\lg \d'\rg)<\frac{3}{26}$.
\item[\rm(4)]  $r\mid q'-1$ and $r\neq 2$:  $\lg z_{q'-1} \rg$ with
$Q(\lg z_{q'-1} \rg)<\frac{1}{26(q'-1)}.$
\item[\rm(5)]  $r\mid q'+1$ and  $r\neq 2,3$: $\lg z_{q'+1}^A (a)\rg$ with  $Q(\lg z_{q'+1}^A(a) \rg)<\frac{1}{13(q'+1)^2};$ and $\lg z_{q'+1}^B \rg$ with
    $Q(\lg z_{q'+1}^B \rg)<\frac{1}{26(q'+1)}.$
\item[\rm(6)] $r\di q'^2-q'+1$ with $r\neq 3$: $\lg z_{q'^2-q'+1} \rg$ and
$Q(\lg z_{q'^2-q'+1} \rg)<\frac{1}{26(q'^2-q'+1)}.$
\end{enumerate}}
\end{lem}
\demo
Suppose $c=3$, $q'\geq 5$ and $q\geq 125$. Then
\vspace{-5pt}$${\small\begin{array}{lll}
\frac{|M|}{|\Omega|}&=&\frac{|M|^2}{|G|}=\frac{2m|\PGU(3,q')|^2}{|\PSU(3,q)|}=\frac{6mq'^6(q'^3+1)^2(q'^2-1)^2}{q^3(q^3+1)(q^2-1)}
\leq\frac{6mq'^{16}(1+q'^{-3})^2}{q^8(1-q^{-2})}<\frac{61m}{10q'^8}.
\end{array}}$$
 By Lemmas~\ref{PGU} and \ref{primeproperties}, we need to discuss the following cases, separately.
\vskip 3mm
(1) Case $r=2$ and $p$ is odd:
By Lemma \ref{primeproperties}, we may set $z=z_2$, while $z_2^f=z_2$.
Thus,
{\small$\Fix(\lg z_2 \rg)=\frac{|\N_G(\lg z_2\rg)|}{|\N_M(\lg z_2\rg)|}=\frac{2m q(q^2-1)(q+1)/3}{q'(q'^2-1)(q'+1)\times 2m}=\frac{q'^2(q'^2+q'+1)(q'^2-q'+1)^2}{3}$,} and
\vspace{-5pt}{\small$$\begin{array}{lll}
Q(\lg z_2 \rg)&=&\frac{|M|}{|\Omega|}\frac{(2-1)\Fix(\lg z_2\rg)}{|\N_{M}(\lg z_2 \rg)|}=\frac{|M|}{|\Omega|}\frac{q'^2(q'^2+q'+1)(q'^2-q'+1)^2}{3\times q'(q'^2-1)(q'+1)\times 2m}\leq\frac{61m}{10q'^8}\frac{q'^4\frac{31}{25}}{6m\frac{24}{25}}
<\frac{15}{11q'^4}<\frac1{26}.
\end{array}$$}
\quad(2) Case $r=p:$
By Lemma \ref{primeproperties}, we may set either $z=z_0$ or $z=z_l'$ where $0\leq l\leq d'-1$.

Suppose that $z=z_0.$ Then
$|\N_{T}(\lg z_0\rg)|=q^3(q+1)(p-1)/3$ and $|\N_{M_1}(\lg z_0\rg)|=q'^3(q'+1)(p-1))$.
This implies that $|\N_{G}(\lg z_0\rg)|\leq2m|\N_{T}(\lg z_0\rg)|=2mq^3(q+1)(p-1)/3$,
and $|\N_{M}(\lg z_0\rg)|\geq e|\N_{M_1}(\lg z_0\rg)|=e q'^3(q'+1)(p-1))$, as $f^{\frac{2m}{e}}\in\C_{G}(M_1)$. Now, we conclude that.
Thus,
\vspace{-5pt}{\small$$\begin{array}{ll}
&Q(\lg z_0 \rg)=\frac{|M|}{|\Omega|}\frac{(p-1)\Fix(\lg z_0\rg)}{|\N_{M}(\lg z_0 \rg)|}=\frac{|M|}{|\Omega|}\frac{(p-1)|\N_G(\lg z_0\rg)|}{|\N_{M}(\lg z_0 \rg)|^2}\leq\frac{|M|}{|\Omega|}\frac{2m(p-1)q^3(q+1)(p-1)}{ 3e^2q'^6(q'+1)^2(p-1)^2} <\frac{61m'^2}{15q'^4}< \frac 1{26}.
\end{array}$$}
\vspace{-5pt}

Suppose that $z=z_l'$ where $0\leq l\leq d'-1$. Since these three elements $z_l'$ are conjugate to each other in $M_1$, we may set $z=z_0'$. Observe that $|\N_{G}(\lg z_0' \rg)|\leq2m|\N_T(\lg z_0\rg)|=2mq^2(p-1)$ and $|\N_{M}(\lg z_0' \rg)|\geq e|\N_{M_1}(\lg z_0'\rg)|=eq'^2(p-1)$. This implies that
\vspace{-5pt}{\small$$\begin{array}{ll}
Q(\lg z_0' \rg)=\frac{|M|}{|\Omega|}\frac{(p-1)\Fix(\lg z_0'\rg)}{|\N_{M}(\lg z_0' \rg)|}=\frac{|M|}{|\Omega|}\frac{(p-1)|\N_{G}(\lg z_0'\rg)|}{|\N_{M}(\lg z_0'\rg)|^2}\leq\frac{|M|}{|\Omega|}\frac{2m(p-1)q'^4}{q'^2(p-1)\times e^2}<\frac{61m}{10q'^8}\frac{2mq'^2}{e^2}=\frac{61m'^2}{5q'^6}<\frac 1{26}.
\end{array}$$}
\quad(3) Case $r=3$ and $r\mid q'+1:$ By Lemma~\ref{PGU}, we take $z=z_A,z_B$ and $\delta'$. Since these are diagonal matrices, we have $\lg z\rg^f=\lg z\rg.$ An it follows from Lemma~\ref{fieldconjugate} that there are two conjugacy classes subgroup in $T$ of order $3$, due to $9\mid q+1$.

$z=z_A:$ $|\N_{G}(\lg z_A\rg)|=2ms(q+1)^2$ and $|\N_{M}(\lg z_A\rg)|=6ms'(q'+1)^2$, where $s,s'\in\{1,2\}$.
Hence,
\vspace{-5pt}{\small$$\begin{array}{lll}
Q(\lg z_A \rg)&=&\frac{|M|}{|\Omega|}\frac{(3-1)\Fix(\lg z_A \rg)}{|\N_{M}(\lg  z_A  \rg)|}=\frac{|M|}{|\Omega|}\frac{2|\N_{G}(\lg z_A\rg)|}{|\N_{M_2}(\lg z_A\rg)|^2}=\frac{|M|}{|\Omega|}\frac{2s(q'^2-q'+1)^2}{3s'\times 3s'(q'+1)^2\times 2m}< \frac{61m}{10q'^8}\frac{2(q'^2-q'+1)^2}{9m(q'+1)^2}\\
&&\leq\frac{61m}{10q'^8}\frac{2q'^2}{9m}=\frac{61}{45q'^6}.
\end{array}$$}
\quad$z=z_B$ or $z=\delta':$
Since $c=3$, we get that $\delta'\in\PSU(3,q)$. Obviously, $\lg \delta'\rg$ is not $M$-conjugate to $\lg z_A\rg$ and $\lg z_B\rg$. But $\lg \delta'\rg$ is $T$-conjugate to $\lg z_B\rg$, where $\delta'\in[w^2,1,1]\Z(\GU(3,q'))$ with
$w\in\FF_{q'^2}^*$ and $|w|=3$.

Since
$|\N_{G}(\lg z_B \rg)|=|\N_{G}(\lg \d'\rg)|=q(q+1)(q^2-1)2m/3$, $|\N_{M}(\lg \d'\rg)|=6m(q'^2-q'+1)$ and $|\N_{M}(\lg z_B\rg)|=2mq'(q'+1)(q'^2-1)$, we have that
\vspace{-5pt}{\small$$
\begin{array}{l}
\begin{array}{ll}
&\Fix(\lg z_B \rg)=\Fix(\lg \d'\rg)=\frac{q(q+1)(q^2-1)2m}{3\times q'(q'+1)(q'^2-1)2m}+\frac{q(q+1)(q^2-1)2m}{(q'^2-q'+1)18m}
\\
&=\frac{1}{3}q'^2(q'^2-1'+1)(q'^2+q'+1) [q'^2-q'+1+\frac{1}{3}q'(q'+1)^2(q'-1)];
\end{array}\\
\begin{array}{lll}
 Q(\lg z_B\rg)&=&\frac{|M|}{|\Omega|}\frac{2\Fix(\lg z_B \rg)}{|\N_{M}(\lg z_B\rg)|}=\frac{|M|}{|\Omega|}\frac{2}{3}\frac{q'^2(q'^2-q'+1)(q'^2+q'+1)[q'^2-q'+1+\frac{1}{3}q'(q'+1)^2(q'-1)]}{q'(q'+1)(q'^2-1)2m}\\
&\leq&\frac{61m}{10q'^8}\frac{q'^3(q'^{-2}+q'^{-1}+1)q'^2[3q'^2-3q'+3+q'(q'+1)^2(q'-1)]}{9mq'^3(1-q'^{-2})}\\
&\leq&\frac{61m}{10q'^8}\frac{\frac{31}{25}q'^2(q'^4+q'^3+2q'^2-4q'+3)}{9m\frac{24}{25}}
\leq\frac{61m}{10q'^8}\frac{899q'^6}{5000m}<\frac{31}{28q'^2};
\end{array}\\
\begin{array}{lll}
 Q(\lg \d'\rg)&=&\frac{|M|}{|\Omega|}\frac{2\Fix(\lg \d'\rg)}{|\N_{M}(\lg \d'\rg)|}=\frac{|M|}{|\Omega|}\frac{2}{3}\frac{q'^2(q'^2-q'+1)(q'^2+q'+1)[q'^2-q'+1+\frac{1}{3}q'(q'+1)^2(q'-1)]}{3(q'^2-q'+1)\times 2m}\\
&\leq&\frac{61m}{10q'^8}\frac{q'^8(q'^{-2}+q'^{-1}+1)(1+q'^{-1}+2q'^{-2}-4q'^{-3}+3q'^{-4})}{27m}
<\frac{2}{5q'^2}.
\end{array}\end{array}$$}

Consequently, we have $Q(\lg z_A\rg)+Q(\lg z_B\rg)+Q(\lg \d'\rg)<\frac{61}{45q'^6}+\frac{31}{28q'^2}+\frac{2}{5q'^2}<\frac{3}{26}.$

\vskip 3mm
\quad(4) Case $r\mid q'-1$ and $r\neq 2:$
By Lemma \ref{primeproperties}, we may set $z=z_{q'-1}$. It follows that
$|\N_{G}(\lg z_{q'-1}\rg)|=4m(q^2-1)/{3}$ and
$|\N_{M}(\lg z_{q'-1}\rg)|=4m(q'^2-1)$.
Therefore, we have
\vspace{-5pt}{\small$$\begin{array}{ll}
&(q'-1)Q(\lg z_{q'-1}\rg)=(q'-1)\frac{|M|}{|\Omega|}\frac{(r-1)\Fix(\lg z_{q'-1}\rg)}{|\N_{M}(\lg z_{q'-1}\rg)|}=(q'-1)\frac{|M|}{|\Omega|}
\frac{(r-1)|\N_{G}(\lg z_{q'-1}\rg)|}{|\N_{M}(\lg z_{q'-1}\rg)|^2}\\
&=(q'-1)\frac{|M|}{|\Omega|}\frac{(r-1)\frac{4m(q^2-1)}{3}}{(4m(q'^2-1))^2}
\leq\frac{61m}{10q'^8}\frac{(r-1)q'^5(q'^{-2}+q'^{-4}+1)}{12m q'^2(1-q'^{-2})}
\leq\frac{61m}{10q'^8}\frac{(r-1)217q'^3}{2400m}
<\frac{4}{7q'^4}\lneqq\frac{1}{26},
\end{array}$$}
\quad(5)
Case $r\mid q'+1$ and $r\neq 2,3:$
Since $r \notin {2, 3}$ and $q' \ge 5$, it follows that $q' \ge 9$.
By Lemma \ref{primeproperties}, we may take $z=z_{q'+1}^A(a)$ and $z=z_{q'+1}^B$, where $2\leq a\leq \frac{r-1}{2}$.
From the proof of Lemma \ref{primeproperties}, we get that these different $M_1$-conjugate classes are also not conjugate in $T$.

$z=z_{q'+1}^A(a):$ It follows from Lemma \ref{primeproperties} that $ |\N_{G}(\lg z_{q'+1}^A(a)\rg)|=\frac{(q+1)^2\times 2m\times s}{3}$ and
$|\N_{M}(\lg z_{q'+1}^A(a)\rg)|$
$= s'(q'+1)^2\times 2m$, where $s,s'$ are divisors of $6$. Hence,
\vspace{-5pt}{\small$$\begin{array}{ll}
&\frac{(q'+1)^2}{2}Q(\lg z_{q'+1}^A(a)\rg)=\frac{(q'+1)^2}{2}\frac{|M|}{|\Omega|} \frac{(r-1)\Fix(\lg z_{q'+1}^A(a)\rg)}{|\N_{M}(\lg z_{q'+1}^A(a)\rg)|}=\frac{(q'+1)^2}{2}\frac{|M|}{|\Omega|} \frac{(r-1)|\N_{G}(\lg z_{q'+1}^A(a)\ rg)|}{|\N_{M}(\lg z_{q'+1}^A(a)\rg)|^2}\\
&=\frac{(q'+1)^2}{2}\frac{|M|}{|\Omega|}\frac{(r-1)s(q'^2-q'+1)^2}{3s'^2(q'+1)^2\times 2m}
\leq \frac{(q'+1)^2}{2}\frac{61m}{10q'^8} \frac{(r-1)(q'^2-q'+1)^2}{m(q'+1)^2}
\leq\frac{61}{20q'^{3}}\lneqq\frac{1}{26}.
\end{array}$$}
\quad$z=z_{q'+1}^B:$ Applying Lemma \ref{primeproperties}, we obtain $|\N_{G}(\lg z_{q'+1}^B\rg)|= \frac{q(q^2-1)(q+1)\times 2m}{3}$ and
$|\N_{M}(\lg z_{q'+1}^B\rg)|=q'(q'^2-1)(q'+1)\times 2m$.
Hence,
\vspace{-5pt}{\small$$\begin{array}{lll}
(q'+1)Q(\lg z_{q'+1}^B\rg)&=&(q'+1)\frac{|M|}{|\Omega|}\frac{(r-1)\Fix(\lg z_{q'+1}^B\rg)}{|\N_{M}(\lg z_{q'+1}^B\rg)|}=(q'+1)\frac{|M|}{|\Omega|}
\frac{(r-1)|\N_{G}(\lg z_{q'+1}^B\rg)|}{|\N_{M}(\lg z_{q'+1}^B\rg)|^2}\\
&=&(q'+1)
\frac{|M|}{|\Omega|}\frac{(r-1)q'^2(q'^2-q'+1)^2(q'^2+q'+1)}{3\times q'(q'^2-1)(q'+1)\times 2m}=\frac{|M|}{|\Omega|}\frac{(r-1)q'(q'^2-q'+1)^2(q'^2+q'+1)}{6m(q'^2-1)}\\
&\leq&\frac{61m}{10q'^8}\frac{(r-1)q'^7(q'^{-2}+q'^{-1}+1)}{6m(1-q'^{-2})q'^2}
\leq\frac{15}{11q'^2}\lneqq\frac{1}{26},
\end{array}$$}
\quad(6) Case $r\mid q'^2-q'+1$ and $r\neq 3:$
Following the same arguments as in the proof of Lemma \ref{QM} for case (6), we can take $z=z_{q'^2-q'+1}$ as a diagonal matrix. Further, by Lemma \ref{primeproperties}, we have that
 $|\N_{G}(\lg z_{q'^2-q'+1}\rg)|=2m(q^2-q+1)$,
$|\N_{M}(\lg z_{q'^2-q'+1}\rg)|=6m(q'^2-q'+1)$.
This follows that
\vspace{-5pt}{\small$$\begin{array}{ll}
&(q'^2-q'+1)Q(\lg z_{q'^2-q'+1}\rg)=(q'^2-q'+1)\frac{|M|}{|\Omega|}\frac{(r-1)(q'^6-q'^3+1)}{6m(q'^2-q'+1)\times(q'^2-q'+1)\times 3}=\frac{|M|}{|\Omega|}\frac{(r-1)(q'^6-q'^3+1)}{18m(q'^2-q'+1)}\\
&\leq\frac{61m}{10q'^8}\frac{(r-1)q'^6}{18mq'^2(1-q'^{-1})}
\leq\frac{61}{144q'^2}\lneqq\frac{1}{26}.\hskip 9.1cm \Box
\end{array}$$}

\begin{lem}\label{3SUF}
Let $c=e=3$, $q=p^m$ and $q'=p^{m'}\geq 5$. Then every cyclic subgroup  $\lg z \rg$ of prime order $r$ of $M$  not contained in $M_1$
is $M$-conjugate to  one of the following groups.
{\small\begin{enumerate}
 \item[\rm(1)]  $r\notin \{e,2\}$: $z=f^{\frac{2m}r}$ with $Q(\lg f^{\frac{2m}r}\rg)< \frac{1}{26m'}.$
\item[\rm(2)] $r=2:$  $z=f^m$ with $Q(\lg f^m \rg)< \frac{1}{26}.$
\item[\rm(3)] $r=e:$ $\lg f^{2m'}\rg$, $z_Af^{2m'}$, $\lg z_Bf^{2m'}\rg$ and $\lg\delta'f^{2m'}\rg$ with
    $Q(\lg f^{2m'} \rg)+Q(\lg z_A f^{2m'}\rg)+Q(\lg z_Bf^{2m'}\rg)+Q(\lg \delta' f^{2m'}\rg)<\frac{1}{26}.$
\end{enumerate}}
\end{lem}
\demo
If $u \in M \setminus M_1$, then we can write $z=u_1f_1$ for some element $u_1\in M_1$ and $f_1\in\lg f \rg$. Since $|z|=r$, we get   $|f_1|=r$ and so we may set $z=u_1f'$, where $f':=f^{\frac{2m}{r}}$. In what follows, we discuss three cases, separately.
\vskip 3mm
(1) Case $r\notin\{2,e\}:$
\vskip 3mm
 Proposition \ref{fieldconjugate} allows us to set $z = f^{\frac{2m}{r}}$, and there is only one $M_1$-class $\lg f^{\frac{2m}{r}}\rg$. We also get that
$|\N_{G}(\lg f^{\frac{2m}{r}}\rg)|=|\PSU(3,q'^{\frac{3}{r}})\rtimes\lg f\rg)|=\frac{2m}{3}(q'^{\frac{3}{r}})^3((q'^{\frac{3}{r}})^3+1)((q'^{\frac{3}{r}})^2-1)$
and
$|\N_{M}(\lg f^{\frac{2m}{r}}\rg)|=|\PGU(3,q'^{\frac{1}{r}})\rtimes\lg f\rg)|=2m(q'^{\frac{1}{r}})^3((q'^{\frac{1}{r}})^3+1)((q'^{\frac{1}{r}})^2-1)$.
Hence,
\vspace{-5pt}{\small$$
\begin{array}{lll}
m'Q(\lg f^{\frac{2m}{r}}\rg)&=&m'\frac{|M|}{|\Omega|}\frac{(r-1)\Fix(\lg f^{\frac{2m}{r}}\rg)}{|\N_{M}(\lg f^{\frac{2m}{r}}\rg)|}=m'\frac{|M|}{|\Omega|}\frac{(r-1)|\N_G(\lg f^{\frac{2m}{r}}\rg)|}{|\N_{M}(\lg f^{\frac{2m}{r}}\rg)|^2}\\
&\leq&\frac{m'|M|}{|\Omega|}\frac{(r-1)q'^{\frac{16}{r}}(q'^{\frac{-9}{r}}+1)(q'^{\frac{-4}{r}}+q'^{\frac{-2}{r}}+1)}{6mq'^{\frac{8}{r}}(1-q'^{\frac{-2}{r}})}
\leq\frac{61m m'}{10q'^{8}}\frac{(r-1)q'^{\frac{8}{5}}(p^{\frac{-9m'}{m'}}+1)(p^{\frac{-4m'}{m'}}+p^{\frac{-2m'}{m'}}+1)}{6m(1-p^{\frac{-2m'}{m'}})}\\
&<&\frac{9m'^2}{5q'^{\frac{32}{5}}}\lneqq\frac{1}{26}.
\end{array}$$}
\quad(2) Case $r=2:$
 \vskip 3mm
 By Proposition \ref{fieldconjugate}, we may take $z=f^{m}$.
 From \cite[Proposition $3.3.15$]{BGbook}, we have
${\small\Fix(\lg f^m\rg)=\frac{|\N_{G}(\lg f^m\rg)|}{|\N_{M}(\lg f^m\rg)|}=\frac{|\PSO(3,q'^3)|2m}{|\PGO(3,q')|2m}=\frac{q'^3(q'^6-1)}{q'(q'^2-1)}=q'^2(q'^4+q'^2+1)}$. Consequently, we conclude that
\vspace{-5pt}{\small$$\begin{array}{lll}
Q(\lg f^m\rg)&=&\frac{|M|}{|\Omega|}\frac{\Fix(\lg f^m\rg)}{|\N_{M_2}(\lg  f^m\rg)|}=\frac{|M|}{|\Omega|}\frac{q'^2(q'^4+q'^2+1)}{q'(q'^2-1)\times 2m}<\frac{4}{q'^5}\lneqq\frac{1}{26}.
\end{array}$$}
\quad(3) Case $r=e=d=d'=3:$
\vskip 3mm
Therefore, $z$ can be expressed as $u_0f^{2m'}$, where $u_0\in M_1$. Since $[x,f^{2m'}]=1$ for any $x\in M_1$, we have $|u_0|=r=3$. Therefore, by Lemma \ref{PGU}, $z$ must be of the form $f^{2m'}$, $z_Af^{2m'}$, $z_Bf^{2m'}$ or $\delta'f^{2m'}$, where $\d'\in M_1\setminus M_0.$

By Proposition \ref{fieldconjugate}, we have that these elements are conjugate to each other in $\PGU(3,q)$.
Now, we will prove that they are also conjugate to each other in $T$.
Since there exists $g\in\PGU(3,q)$ such that $f'^g=g^{-1}g^{f'^{-1}}f'=z_Af'$. It follows that $g=g^{f'^{-1}}(z_A)^{-1}$. Followed from \cite[Lemma 2.1.1]{BGbook} that $x\in\PSU(3,q)$ if and only if $\det(x)=\lambda^3$ for some $\lambda\in\FF_{q^2}^*$. Suppose that $\det(g)=\xi^k$, where $\FF_{q^2}^*=\lg \xi\rg$. Now, we have $\xi^k=\xi^{-kq'^2}\xi^{3l}$, where $k,l$ are integers, which implies that $3\mid k(1+q'^2)$. It follows that $3\mid k$, as $(q'+1,3)=3$ and $(q'^2+1,3)=1$. Therefore, $g\in\PSU(3,q)$. Finally, by analogous arguments, the elements $f'$, $z_A f'$, $z_B f'$, and $\delta' f'$ are all conjugate to each other within $\mathrm{PSU}(3,q)$.
Therefore, by Proposition \ref{man}, we have
 \vspace{-5pt}$${\small\begin{array}{ll}
&|\Fix(\lg f'\rg)|=|\Fix(\lg z_A f'\rg)|=|\Fix(\lg z_B f'\rg)|=|\Fix(\delta'f')|\\
&=\frac{|\N_{G}(\lg f'\rg)|}{|\N_{M}(\lg f'\rg)|}+\frac{|\N_{G}(\lg z_Af'\rg)|}{|\N_{M}(\lg z_Af'\rg)|}+\frac{|\N_{G}(\lg z_Bf'\rg)|}{|\N_{M}(\lg z_Bf'\rg)|}+\frac{|\N_{G}(\lg \delta'f'\rg)|}{|\N_{M}(\lg \delta'f'\rg)|}\\
&=\frac{|\N_{G}(\lg f'\rg)|}{|\N_{M}(\lg f'\rg)|}
+\frac{|\N_{G}(\lg f'\rg)|}{|\C_{M}(\lg z_A\rg)|}
+\frac{|\N_{G}(\lg f'\rg)|}{|\C_{M}(\lg z_B\rg)|}
+\frac{|\N_{G}(\lg f'\rg)|}{|\C_{M}(\lg \delta'\rg)|}
\\
&=\frac{|\PGU(3,q'):\lg f\rg|}{|\PGU(3,q'):\lg f\rg|}
+\frac{|\PGU(3,q'):\lg f\rg|}{[3(q'+1)^2]r}
+\frac{|\PGU(3,q'):\lg f\rg|}{[q'(q'+1)(q'^2-1)]r}
+\frac{|\PGU(3,q'):\lg f\rg|}{[3(q'^2-q'+1)]r}\\
&=1+\frac{q'^3(q'^3+1)(q'^2-1)\times 2m}{(q'+1)^2\times 3r}+\frac{q'^3(q'^3+1)(q'^2-1)\times 2m}{q'(q'+1)(q'^2-1)r}+\frac{q'^3(q'^3+1)(q'^2-1)\times 2m}{3(q'^2-q'+1)r}\\
&=1+2m'q'^2(\frac{2}{3}q'^4-\frac{1}{3}q'^3+\frac{4}{3}q'^2-\frac{5}{3}q'+1)\leq\frac{4m'q'^6}{3}.
\end{array}}$$
We therefore have
\vspace{-5pt}{\small$$\begin{array}{ll}
&\begin{array}{lll}
Q(\lg f'\rg)&=& \frac{|M|}{|\Omega|}\frac{2\Fix(\lg f'\rg)}{|\N_{M}(\lg f'\rg)|}\leq \frac{|M|}{|\Omega|}\frac{2\times\frac{4m'q'^6}{3}}{q'^3(q'^3+1)(q'^2-1)\times 2m}\leq\frac{|M|}{|\Omega|}\frac{25}{54q'^2},
\end{array}$$\\
&\begin{array}{lll}
Q(\lg z_Af'\rg)&=&\frac{|M|}{|\Omega|}\frac{2\Fix(\lg z_Af'\rg)}{|\N_{M}(\lg z_Af'\rg)|}\leq\frac{|M|}{|\Omega|} \frac{2\times\frac{4m'q'^6}{3}}{(q'+1)^2\times 9}\leq\frac{|M|}{|\Omega|}\frac{8m'q'^4}{27},\\
\end{array}\\
&\begin{array}{lll}
Q(\lg z_Bf'\rg)&=&\frac{|M|}{|\Omega|}\frac{2\Fix(\lg z_Bf'\rg)}{|\N_{M}(\lg z_Bf'\rg)|}\leq\frac{|M|}{|\Omega|} \frac{2\times\frac{4m'q'^6}{3}}{3q'(q'+1)(q'^2-1)}\leq\frac{|M|}{|\Omega|}\frac{25m'q'^2}{27},\\
\end{array}\\
&\begin{array}{lll}
Q(\lg \delta'f'\rg)&=&\frac{|M|}{|\Omega|}\frac{2\Fix(\lg \delta'f'\rg)}{|\N_{M}(\lg \delta'f'\rg)|}\leq\frac{|M|}{|\Omega|} \frac{2\times\frac{4m'q'^6}{3}}{9(q'^2-q'+1)}\leq\frac{|M|}{|\Omega|}\frac{10m'q'^4}{27},\\
\end{array}\\
&\begin{array}{ll}
&Q(\lg f' \rg)+Q(\lg z_A f'\rg)+Q(\lg z_Bf'\rg)+Q(\lg \delta' f'\rg)<\frac{61m}{10q'^8}(\frac{25}{54q'^2}+\frac{18m'q'^4}{27}+\frac{25m'q'^2}{27})\\
&\leq \frac{305m'}{36q'^{10}}+\frac{61m'^2}{5q'^4}+\frac{305m'^2}{18q'^6}\lneqq\frac{1}{26}.
\end{array}
\end{array}$$}

\qed

\begin{theorem}\label{c=3proof}
Let $c=e=3$, $G=\PSigmaU(3,q)$ and $M=\PGU(3,q'):\lg f\rg$. Then the BG-Conjecture holds for $G$ acting on $[G:M]$ by right multiplication.
\end{theorem}
\demo
Taking into account, all conjugate classes of subgroups $\lg z\rg $  of prime order in $M$ and an upper bound for the number $|Q(\lg z\rg )|$ have been determined  in  Lemmas~\ref{3SU} and \ref{3SUF},
  we are ready to show $\check{Q}(G)< \frac 12$ so that  the BG-Conjecture holds.
\vskip 3mm

Write $m'=r_1^{k_1} r_2^{k_2}\cdots r_{l}^{k_l}$ if $m'\neq 1$, where $2\le r_1\lneqq r_2\lneqq \cdots \lneqq r_l$ are distinct primes and $k_i\geq 1$ for all $i$. Set $m_i=\frac{m}{r_i}$, $q_i=p^{m_i}$ and $f_i=f^{m_i}.$
Further, let
\vspace{-5pt}\begin{enumerate}
  \item  [{\rm(a)}] $q'-1=2^{\epsilon_2}t_1^{j_1} t_2^{j_2}\cdots t_{o}^{j_o}$, where $3\le t_1\lneqq t_2\lneqq \cdots \lneqq t_o$ are distinct primes;
    \item  [{\rm(b)}] $q'+1=2^{\epsilon_1}3^{\kappa_1}s_1^{i_1} s_2^{i_2}\cdots s_{h}^{i_h}$,
 where $5\le s_1\lneqq s_1\lneqq \cdots \lneqq s_h$ are distinct primes; and
  \item [{\rm(c)}] $q'^2-q'+1=2^{\epsilon_3}3^{\kappa_3}u_1^{\ell_1} u_2^{\ell_2}\cdots u_{\imath}^{\ell_{\imath}}$ where $5\le \ell_1\lneqq \ell_2\lneqq \cdots \lneqq \ell_\imath$ are distinct primes.
\end{enumerate}
Furthermore, denote by $z_{q'-1, i}$ an element of order $t_i$ in $\lg z_{q'-1}\rg$; denote by $z_{q'+1, i}^A(a)$ and $z_{q'+1, i}^B$ elements of order $s_i$ in $\lg z_{q'+1}^A(a)\rg$ and $\lg z_{q'+1}^B\rg$, respectively; denote
 by $z_{q'^2-q'+1, i}$ an element of order $u_i$ in $\lg z_{q'^2-q'+1}\rg$.
Recall $\mathcal{P}^*(M):=\{\lg g_1\rg,\lg g_2\rg,\cdots,\lg g_k\rg\}$ is the set of representatives  of conjugacy classes of subgroups of prime order in $M$.

By Lemmas \ref{3SU} and \ref{3SUF}, we have
\vspace{-5pt}$${\small\begin{array}{ll}
 &\check{Q}(G)=\frac{|M|}{|\O|}(\sum_{i=1}^k(|\lg g_i\rg|-1)\frac{|\Fix(\lg g_i\rg)|}{|\N_M(\lg g_i\rg )|})=\sum_{i=1}^kQ(\lg g_i\rg)\\
 &\leq  (Q(\lg z_0 \rg)+Q(\lg z_0'\rg)
+Q(\lg z_2 \rg)+Q(\lg z_3^A\rg)
+Q(\lg z_3^B \rg)+Q(\lg \delta'\rg))\\
&+\sum\limits_{i=1}^{o}Q(\lg z_{q'-1, i}\rg)+
\sum\limits_{i=1}^{h}\sum\limits_{1\le a\le \frac{s_i-3}2} Q(\lg z_{q'+1, i}^A(a)\rg)+\sum\limits_{i=1}^{h}Q(\lg z_{q'+1, i}^B\rg)
+\sum\limits_{i=1}^{\imath} Q(\lg z_{q'^2-q'+1, i}\rg)\\
&+\sum\limits_{r_i\mid m,r_i\neq 2, e}Q(\lg f^{\frac{2m}{r_i}} \rg)
+(Q(\lg f^{2m'} \rg)+Q(\lg z_A f^{2m'}\rg)+Q(\lg z_Bf^{2m'}\rg)+Q(\lg \delta' f^{2m'}\rg))+Q(\lg f^m\rg)\\
&\le  \frac{6}{26}+(q'-1) Q(\lg z_{q'-1}\rg)+ \frac{(q'+1)^2}{2} Q(\lg z_{q'+1}^A(a)\rg)+(q'+1)Q(\lg z_{q'+1}^B\rg)\\
&+(q'^2-q'+1)Q(\lg z_{q'^2-q'+1}\rg)+m'Q(\lg f^{\frac{2m}{r}}\rg)+\frac 1{26}+\frac 1{26}  \\
&\le  \frac{6}{26}+7\cdot\frac{1}{26}<\frac 12,
\end{array}}$$
as desired.$\hskip 13.8cm \Box$

\section{$M_0\cong\PSL(2,7)$ or $\PSL(2,9)$}
\begin{theorem}\label{psl27}
 Let $T=\PSU(3,q)$ and $M_0=\PSL(2,7)$ with $q=p\equiv 3,5,6 \pmod 7$ and $q\neq 5$, so that  $q=p\geq 13$. Let $T\leq G\leq \Aut(\PSU(3,q))$ and $M$  a maximal subgroup of $G$ containing $M_0$.
Then the BG-Conjecture holds for $G$ acting on $[G:M]$ by right multiplication.
\end{theorem}
\demo  By checking  \cite[Tabel 8.6]{BGbook}, we may set  either $G=T$ and $M=M_0$; or  $G=\PSigmaU(3,q)$  and $M=M_0:\lg f \rg\cong \PGL(2,7)$.  By Remark~\ref{big},  it suffices to show $ \check{Q}(G)<\frac 12$  for the second case.

Set $\mathcal{P}^*(M)$ be the set of all the representatives of conjugacy classes of subgroups of prime order in $M$.
 Then \vspace{-5pt}$$\mathcal{P}^*(M)=\{\lg x_{1}\rg, \lg x_{2}\rg,  \lg x_3\rg, \lg x_4\rg\},$$
 where $x_1, x_3, x_4\in M_0$ and $x_2\in M\setminus M_0$,
  $|x_1|=2$,  $|x_{2}|=2$,  $|x_3|=3$ and $|x_4|=7$. Note that $\lg x_1\rg$ and $\lg x_2\rg$ are also not conjugate to each other in $G$. Moreover, by Lemma \ref{primeproperties} we have $\Fix(\lg x_i\rg)=\frac{|\N_G(\lg x_i\rg)|}{|\N_{M}(\lg x_i\rg)|}$, that is
  \vspace{-5pt}$${\small\begin{array}{llllll} \Fix(\lg x_1\rg)&=&\frac{2\times q(q^2-1)(q+1)/d}{16}=\frac{q(q^2-1)(q+1)}{8d};
&\Fix(\lg x_2\rg)&=&\frac{2|\SO(3,q)|}{12}=\frac{q(q^2-1)}{6};\\
  \Fix(\lg x_3\rg)&\le& \frac{2\times q(q+1)(q^2-1)/d}{12}\leq\frac{7q^4}{39d};
 &\Fix(\lg x_4\rg)&\leq & \frac{2\times q(q+1)(q^2-1)/d}{42}\leq\frac{2q^4}{39d}.\\
 \end{array}}$$
 Therefore, \vspace{-5pt}{\small$$\begin{array}{lll} \check{Q}(G)&:=&\frac{|M|}{|\O|}(\sum_{i=1}^4(|\lg x_i\rg|-1)\frac{\Fix(\lg x_i\rg)}{|\N_M(\lg x_i\rg )|})\\
&\le &\frac{336^2\times d}{2\times q^3(q^3+1)(q^2-1)}(\frac{q(q^2-1)(q+1)}{8d\times 16}+\frac{q(q^2-1)}{6\times 12}+\frac{2\times7q^4}{39d\times12}+\frac{6\times2q^4}{39d\times42})<\frac{1}{2} . \hskip 3cm \Box\\
 \end{array}$$}

\begin{theorem}\label{psl29}
Let $T=\PSU(3,q)$ and $M_0=\PSL(2,9)$ with $q=p\equiv 11,14 \pmod{15}$.   Let $T\leq G\leq \Aut(\PSU(3,q))$ and $M$ is  a maximal subgroup of $G$ containing $M_0$.
Then the BG-Conjecture holds for $G$ acting on $[G:M]$.
\end{theorem}
\demo   By checking  \cite[Tabel 8.6]{BGbook} where $q\ge 7$, we get that either $G=T$ and $M=M_0$; or  $G=\PSigmaU(3,q)\cong T:\lg f\rg$  and $M=M_0:\lg f \rg$, where  $f$ is given  by
  $g^f=(g^{-1})^T$ for any $g\in T$  so that  $M\cong \PGL(2,9)$ (refer to \cite[P.34]{Bray}). By Remark~\ref{big},  it suffices to show $ \check{Q}(G)<\frac 12$  for the second case.

   Set $\mathcal{P}^*(M)$ be the set of all the representatives of conjugacy classes of subgroups of prime order in $M$.
 Thus, \vspace{-5pt}$$\mathcal{P}^*(M)=\{\lg x_{1}\rg, \lg x_{2}\rg,  \lg x_3\rg, \lg x_4\rg\},$$ where
 where $x_1, x_3, x_4\in M_0$ and $x_2\in M\setminus M_0$,
  $|x_1|=2$,  $|x_{2}|=2$,  $|x_3|=3$ and $|x_4|=5$.  Note that $\lg x_1\rg$ and $\lg x_2\rg$ are also not conjugate to each other in $G$. Moreover, by Lemma \ref{primeproperties} and the fact $q\geq11$, we have
\vspace{-5pt}$${\small\begin{array}{llllll} \Fix(\lg x_1\rg)&=&\frac{2\times q(q^2-1)(q+1)/d}{16}=\frac{q(q^2-1)(q+1)}{8d};
&\Fix(\lg x_2\rg)&=&\frac{2|\SO(3,q)|}{20}=\frac{q(q^2-1)}{10};\\
  \Fix(\lg x_3\rg)&\leq &\frac{2\times q(q+1)(q^2-1)/d}{18}\leq\frac{4q^4}{33d};
 &\Fix(\lg x_4\rg)&\leq &\frac{2\times q(q+1)(q^2-1)/d}{20}\leq\frac{6q^4}{55d}.\\
 \end{array}}$$
 Therefore, \vspace{-5pt}{\small$$\begin{array}{lll} \check{Q}(G):&=&\frac{|M|}{|\O|}(\sum_{i=1}^4(|\lg x_i\rg|-1)\frac{\Fix(\lg x_i\rg)}{|\N_M(\lg x_i\rg )|})\\
&\le& \frac{336^2\times d}{2\times q^3(q^3+1)(q^2-1)}(\frac{q(q^2-1)(q+1)}{8d\times 16}+\frac{q(q^2-1)}{10\times 20}+\frac{2\times 4q^4}{33d\times 18}+\frac{4\times 6q^4}{55d\times 20})<\frac 12.\hskip 3cm \Box \\
 \end{array}$$}

\vskip 3mm
\begin{center}{\large\bf Acknowledgements}\end{center}
\vskip 2mm
The first author thanks the supports of the National Natural Science Foundation of China (12301446, 12571362).
The second author thanks the supports of the National Natural Science Foundation of China (12471332).
The third author thanks the supports of the National Natural Science Foundation of China (12301422).

\end{document}